\documentclass[reqno,10pt,a4paper]{amsart}

\usepackage{amssymb}
\usepackage{tikz,tikz-3dplot,subcaption}
\usepackage{tikz-cd} 
\usepackage[colorlinks=false]{hyperref}
\usepackage[abbrev,bibtex-style]{amsrefs}

\title[A conformal Hopf--Rinow theorem]{A conformal Hopf--Rinow theorem for semi-Riemannian spacetimes}

\author[A.~Burtscher]{Annegret Burtscher}

\thanks{Department of Mathematics, IMAPP, Radboud University, PO Box 9010, 6500 GL Nijmegen, The Netherlands, \textit{Email}: \texttt{burtscher@math.ru.nl}.\\
The final version of this article is published open access in \emph{Indagationes Mathematicae} and available at \url{https://doi.org/10.1016/j.indag.2025.01.004}.}



\date{\today}

\numberwithin{equation}{section}

\theoremstyle{plain}
\newtheorem{thm}{Theorem}[section]
\newtheorem{cor}[thm]{Corollary}
\newtheorem{prop}[thm]{Proposition}
\newtheorem{lem}[thm]{Lemma}

\theoremstyle{definition}
\newtheorem{defn}[thm]{Definition}
\newtheorem{ex}[thm]{Example}

\theoremstyle{remark}
\newtheorem{rem}[thm]{Remark}

\def\vf{\mathfrak{X}}

\def\RR{\mathbb{R}}
\def\NN{\mathbb{N}}
\def\CC{\mathbb{C}}
\def\A{\hat{\mathcal{A}}}

\def\inte{\operatorname{int}}
\def\spann{\operatorname{span}}

\tikzset{
dot/.style = {circle, fill, minimum size=#1,
              inner sep=0pt, outer sep=0pt},
dot/.default = 4pt 
}

\tdplotsetmaincoords{70}{45}
\tdplotsetrotatedcoords{-90}{180}{-90}

\tikzset{surface/.style={draw=blue!70!black, fill=blue!40!white, fill opacity=1}}

\newcommand{\coneback}[4][]{
  \draw[canvas is xy plane at z=#2, #1] (45-#4:#3) arc (45-#4:225+#4:#3) -- (O) --cycle;
  }
\newcommand{\conefront}[4][]{
  \draw[canvas is xy plane at z=#2, #1] (45-#4:#3) arc (45-#4:-135+#4:#3) -- (O) --cycle;
  }

\begin{document}


\begin{abstract}
 The famous Hopf--Rinow Theorem states, amongst others, that a Riemannian manifold is metrically complete if and only if it is geodesically complete. The Clifton--Pohl torus fails to be geodesically complete proving that this theorem cannot be generalized to compact Lorentzian manifolds. On the other hand, Hopf and Rinow characterized metric completeness also by properness. Garc\'ia-Heveling and the author recently obtained a Lorentzian completeness-compactness result for open manifolds with a similar flavor. In this manuscript, we extend the null distance used in this approach and our theorem to proper cone structures and to a new class of semi-Riemannian manifolds, dubbed $(n-\nu,\nu)$-spacetimes. Moreover, we demonstrate that our result implies, and hence generalizes, the metric part of the Hopf--Rinow Theorem.
\end{abstract}

\maketitle

\section{Introduction}
\label{intro}

 In 1931 H.\ Hopf and his student W.\ Rinow published their famous completeness result for Riemannian manifolds. Throughout we consider a manifold to be smooth, paracompact, Hausdorff, connected, finite-dimensional, and without boundary.

\begin{thm}[Hopf--Rinow~\cite{HR}]\label{HRthm}
 Let $(\Sigma,\sigma)$ be a Riemannian manifold. The following statements are equivalent:
 \begin{enumerate}
  \item[(a$_0$)] $(\Sigma,\sigma)$ is geodesically complete for some $p \in \Sigma$, i.e., all geodesics through $p$ are defined for all times.
  \item[(a)] $(\Sigma,\sigma)$ is \emph{geodesically complete}, i.e., all geodesics are defined for all times. 
  \item[(b)] $(\Sigma,d_\sigma)$ is a \emph{complete metric} space.
  \item[(c)] $(\Sigma,d_\sigma)$ satisfies the \emph{Heine--Borel property} (also called \emph{proper} or \emph{boundedly compact}), i.e., every closed and bounded subset of $\Sigma$ is compact.
 \end{enumerate}
\end{thm}

 Originally, Hopf and Rinow called (a$_0$) the abatability postulate (\emph{Abtragbarkeitspostulat}), (a) the infinity postulates (\emph{Unendlichkeitspostulat}), (b) the completeness postulate (\emph{Vollst\"andigkeitspostulat}), and (c) the compactness postulate (\emph{Kompaktheitspostulat}). Most differential geometers focus on and apply the equivalence (a)$\Longleftrightarrow$(b). The metric equivalence (b)$\Longleftrightarrow$(c) already became more central in a generalization of the Hopf--Rinow Theorem to locally compact length metric spaces by S.\ Cohn-Vossen. We refer to Theorem~\ref{HRthm} (b)$\Longleftrightarrow$(c) also as the \emph{metric Hopf--Rinow Theorem} or the \emph{Riemannian completeness-compactness Theorem}. 
 
 For semi-Riemannian manifolds geodesic completeness is rare and difficult to obtain. Even compact Lorentzian manifolds need not be geodesically complete (the Clifton--Pohl torus is a counterexample, see \cite{ON}*{p.\ 193} or \cite{GHL:Rie}*{Ex.\ 2.141}), and strong additional conditions such as homogeneity (obtained by Marsden \cite{Ma72} in 1972) or flatness (obtained by Carri\`ere \cite{Ca89} in 1989) are required to imply it in general settings. The compact case shows that (c)$\Longrightarrow$(a) is false in semi-Riemannian geometry. Since no canonical distance function exists on semi-Riemannian manifolds, a generalization of the Hopf--Rinow Theorem involving (b) was also considered to be impossible in this setting.
 
 The aim of the present manuscript is to show that (and in which sense precisely) the equivalence (b)$\Longleftrightarrow$(c) is robust when leaving the positive definite Riemannian world in negative directions.
 
 We proceed in three steps.
 First, we prove that our recent Lorentzian result \cite{BGH2}*{Theorem 1.4]} is a generalization of Theorem~\ref{HRthm} (b)$\Longleftrightarrow$(c) in Section \ref{sec:LRequiv} (the necessary Lorentzian background is provided in Section \ref{sec:Lbasics}). Second, we extend the relevant definitions and this result to proper cone structures in Section \ref{sec:cone}. Third, we introduce the notion of semi-Riemannian spacetimes, discuss their existence and other properties, and show that they carry proper cone structures in Section \ref{sec:sR}. All steps combined imply, in particular, that the metric Hopf--Rinow Theorem~\ref{HRthm} (b)$\Longleftrightarrow$(c) generalizes to a large class of (noncompact) semi-Riemannian manifolds.
 In the remaining part of this introduction we present and discuss our  results in more detail. 

 \medskip\noindent
 {\bf The Lorentzian setting.} In recent work L.\ Garc\'ia-Heveling and the author proved the following Lorentzian completeness-compactness theorem with \emph{one} negative tangent direction (see Section \ref{sec:Lbasics} for all notions). 

\begin{thm}[Burtscher--Garc\'ia-Heveling \cite{BGH2}]\label{BGHthm}
 Let $(M,g)$ be \emph{spacetime}, i.e., a time oriented Lorentzian manifold. The following are equivalent:
 \begin{enumerate}
  \item[(b')] There exists a time function $\tau \colon M \to \RR$ such that $M$ equipped with the corresponding null distance $\hat d_\tau$ is a complete metric space.
  \item[(c')] $(M,g)$ is globally hyperbolic.
 \end{enumerate}
\end{thm}

 Can we view Theorem~\ref{BGHthm} as the Lorentzian analogue or generalization of the metric Hopf--Rinow Theorem~\ref{HRthm}?
 
 In the present manuscript we show that the answer is yes. At first, the metric Hopf--Rinow Theorem immediately appears to be more special and fragile than Theorem~\ref{BGHthm}. This is because the notions (time functions, null distance, global hyperbolicity) relevant for (b') and (c') are \emph{conformal}ly invariant while the Riemannian distance function is only preserved by isometries. In Section~\ref{sec:LRequiv} we demonstrate that Theorem~\ref{BGHthm} is indeed the correct extension of the metric Hopf--Rinow Theorem~\ref{HRthm} (b)$\Longleftrightarrow$(c) by proving that both properties are (independently) equivalent to their Lorentzian counterparts when restricted to Lorentzian products.  

\begin{thm}\label{prop:equivLR}
 Let $(\Sigma,\sigma)$ be a Riemannian manifold and $(\Sigma',\sigma')=(\RR\times \Sigma, -dt^2 \oplus \sigma)$ be the corresponding Lorentzian product. Then $(\Sigma',\sigma')$ is a stably causal spacetime and the following statements hold:
 \begin{enumerate}
  \item $(\Sigma,d_\sigma)$ satisfies the Heine--Borel property if and only if the causal diamonds of $(\Sigma',\sigma')$ are compact, i.e., $(\Sigma',\sigma')$ is globally hyperbolic.
  \item $(\Sigma,d_\sigma)$ is a complete metric space if and only if $(\Sigma',\hat d_t)$ is a complete metric space with respect to the canonical time function $t(p_0,p_\Sigma) = p_0$.
 \end{enumerate}
 \end{thm}

 What about geodesic completeness?
 In the product case it is known that (a) geodesic completeness of $(\Sigma,\sigma)$ is equivalent to (a') geodesic completeness of $(\Sigma',\sigma')$, and also equivalent to global hyperbolicity of $(\Sigma',\sigma')$ (see, for instance, \cite{BEE}*{Theorems 3.66 and 3.67}). Thus in the Lorentzian product case (and only then!) all properties are equivalent independent of each other:
 
 
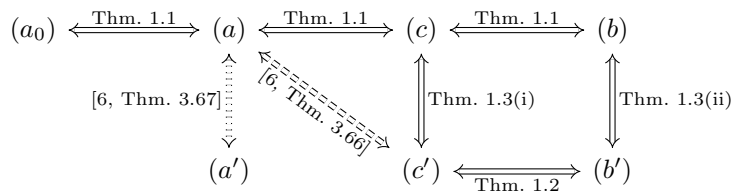
\[ \begin{tikzcd}[row sep=huge, column sep = huge]
(a_0) \arrow[Leftrightarrow]{r}{\text{Thm.~\ref{HRthm}}} &(a) \arrow[Leftrightarrow]{r}{\text{Thm.~\ref{HRthm}}} \arrow[Leftrightarrow,dotted]{d}[swap]{\text{\cite{BEE}*{Thm.~3.67}}} \arrow[Leftrightarrow,dashed]{rd}[sloped,below]{\text{\cite{BEE}*{Thm.~3.66}}} & (c) \arrow[Leftrightarrow]{d}{\text{Thm.~\ref{prop:equivLR} (i)}} \arrow[Leftrightarrow]{r}{\text{Thm.~\ref{HRthm}}} & (b) \arrow[Leftrightarrow]{d}{\text{Thm.~\ref{prop:equivLR} (ii)}} \\%
& (a') & (c') \arrow[Leftrightarrow]{r}[swap]{\text{Thm.~\ref{BGHthm}}} & (b')
\end{tikzcd}
\]
 
 While the equivalence of (a)$\Longleftrightarrow$(c') carries over to warped products $M \times_f \Sigma$ for $(M,g)$ a spacetime, even for warped products with compact Riemannian slices (a)$\not\Longrightarrow$(a') (see \cite{BEE}*{Theorem 3.68} and preceding discussion). More generally, the singularity theorems of Penrose and Hawking (see, for instance, \cite{HE}*{Chapter 9}) from the 1960s reveal that (c')$\not\Longrightarrow$(a') and due to Anti-de Sitter space (a')$\not\Longrightarrow$(c'). It was very recently shown by S\'anchez \cite{San:gh}*{Section 6.4} that for general spacetimes also (c')$\not\Longrightarrow$(a) for its slices. These observations show that on semi-Riemannian manifolds we are forced to let go of geodesic completeness (a) and (a') also in the noncompact setting, and that new tools are needed.
 
  \medskip\noindent
 {\bf The semi-Riemannian setting.}
 That lifting and projecting the properties relevant for the metric Hopf--Rinow Theorem~\ref{HRthm} is doable in the Lorentzian setting raises the hope that it is possible to iterate this procedure and add \emph{arbitrarily many} negative tangent directions. The first obstacle we face is the lack of an analogue of causality theory for semi-Riemannian manifolds, meaning that none of the notions used in (b') and (c') yet exist. To this end we introduce and analyze in Section~\ref{sec:sR} a new subclass of semi-Riemannian manifolds by generalizing the concept of a Lorentzian $(n-1,1)$-spacetime to semi-Riemannian manifolds of any index $0 \leq \nu \leq n$ as follows.

\begin{defn}\label{def:nnuspacetime}
 Let $(M,g)$ be a semi-Riemannian manifold with constant index $0 \leq \nu \leq n = 
 \dim M$. We say that $M$ is \emph{time frame 
 orientable} if it admits $\nu$ continuous vector fields $X_i \in \vf(M)$ that in each $p \in M$ span a negative definite $\nu$-dimensional subspace of $T_pM$.
 
 If $(M,g)$ is time frame orientable and equipped with a fixed set $X = \{X_i \, ; \, i=1,\ldots,\nu\}$ of such vector fields, we say that it is \emph{time frame oriented} and call $(M,g,X)$ a \emph{semi-Riemannian spacetime} or, more specifically, a \emph{$(n-\nu,\nu)$-spacetime}.
\end{defn}
 
 Whether a given manifold $M$ admits such a spacetime structure is a purely topological question. In 1935 E.\ Stiefel asked in his seminal thesis \cite{Sti:parallel}, supervised by H.\ Hopf, precisely this question in the closed manifold case (translated from German):
 \begin{quote}
  Is there a system of $\nu$ continuous vector fields on $M$, 
  which are linearly independent in each point of $M$?
 \end{quote} 
 Such a system of of vector fields is called a \emph{tangent $\nu$-frame} in the literature, and its existence equivalent to the existence of a semi-Riemannian $(n-\nu,\nu)$-spacetime structure on $M$. (This is also the reason why we call it a time frame in Definition \ref{def:nnuspacetime}, even though the order of the $X_i$ is irrelevant for us.)
  
 For $\nu=1$, the Poincar\'e--Hopf Theorem \cite{Ho:vf} from 1927 fully answers the above question in the closed manifold case, characterizing existence by $\chi(M)=0$. In 1955 Markus \cite{Mar} showed that no conditions are needed in the open case.
 
 For $\nu>1$ Stiefel's question is in its generality still open, but led to many exciting developments in algebraic topology. The first necessary conditions were formulated in terms of Stiefel--Whitney classes \cites{Sti:parallel,Whi:top}, and later characteristic classes. In the 1960s and 1970s several breakthrough results were obtained for closed manifolds in special cases (in particular, for $\nu=2,3$, such as the work of Atiyah and Dupot \cites{AtDu,Du:K}, and $\nu=n$), which we discuss in more detail in Section \ref{ssec:sRspacetimes}. 
 More recently, B\"okstedt--Dupont--Svane \cite{BDS:cobord} also give a brief overview and present new results for $\nu=4,5,6,7$ (under additional assumptions on $M$ and $n$). 

 Note that the above investigations all concern the closed case and often allow finite singularities. From our semi-Riemannian perspective, we also immediately understand that the situation for $\nu>1$ is more subtle than for $\nu=1$ because \emph{not} every manifold admitting a semi-Riemannian metric also admits a (possibly different) semi-Riemannian spacetime metric ($\mathbb{S}^2 \times \mathbb{S}^2$ is a counterexample, see Example~\ref{ex:s2}). Nonetheless, it is this open case that we are interested in and ultimately even forced to consider in the present paper, simply because closed manifolds do not admit time functions. Thus the tools developed in this paper may in the future also help to shed light on the existence problem of a tangent $\nu$-frame in the noncompact case. By establishing the following result we show, in particular, that $(n-\nu,\nu)$-spacetimes can be studied with techniques from the theory of cone structures developed in recent years by Fathi and Siconolfi \cites{FaSi,Fa:time}, Bernard and Suhr \cites{BeSu:time,BeSu:gh}, and Minguzzi \cites{Min:cone}.

\begin{thm}\label{thm:sRcones}
 A $(n-\nu,\nu)$-spacetime with index $0<\nu \leq n$ admits a conformally invariant continuous proper cone structure $(M,C)$.
\end{thm}
 
 In the framework of cone structures the notions of time functions and global hyperbolicity appearing in Theorem~\ref{BGHthm} have already been introduced and are well studied. In Section~\ref{sec:cone} we show that also the null distance can be extended to proper cone structures.
 The main result of this manuscript is an extension of Theorem~\ref{BGHthm}, and thus of the metric Hopf--Rinow Theorem, to cone structures.

\begin{thm}\label{conethm}
 Let $(M,C)$ be a proper cone structure. The following are equivalent:
 \begin{enumerate}
  \item[(b'')] There exists a time function $\tau \colon M \to \RR$ such that $(M,\hat d_\tau)$ is a complete metric space.
  \item[(c'')] $(M,C)$ is globally hyperbolic.
 \end{enumerate}
\end{thm}

 Together with Theorem~\ref{thm:sRcones} this characterization yields a conformal com\-plete\-ness-compactness theorem in semi-Riemannian geometry.

\begin{cor}\label{semiRcor}
 For any $(n-\nu,\nu)$-spacetime with $0<\nu \leq n$ we have that (b'') is equivalent to (c'').
\end{cor}

 Above we have already seen that the Riemannian equivalence (b)$\Longleftrightarrow$(c) is a special case of the Lorentzian equivalence (b')$\Longleftrightarrow$(c'). Can we in the same way relate $\nu$ and $\nu+1$ if $\nu \geq 1$?
 
 In Section~\ref{ssec:sRGH} we show that the answer is no. For one, if we extend a $(n-\nu,\nu)$-spacetime $(M,g)$ orthogonally to a $(n-\nu,\nu+1)$-spacetime $(M',g')=(\RR\times M,-dt^2\oplus g)$, then the corresponding null distances cannot be directly related because they are \emph{both} conformally invariant. In fact, even the existence of a null distance on $(M',g')$ requires \emph{additional} properties and is not a given as in the Lorentzian product case. For the same reason global hyperbolicity does not carry over. Nonetheless, we can still recover parts of Theorem~\ref{prop:equivLR} in semi-Riemannian geometry as summarized in the following result obtained in Section~\ref{sec:sR}.

\begin{thm}\label{thm:equivsR}
 Let $(M,g)$ be a $(n-\nu,\nu)$-spacetime, $0 < \nu \leq n$, with time frame orientation defining vector fields $X_1,\ldots,X_\nu$. Then one can orthogonally extend the metric and time frame orientation (with $X_{\nu+1} = \partial_t$) on $M' = \RR \times M$ such that the following implications hold:
 \begin{align*}
  &(M',-dt \oplus g) \text{ is a globally hyperbolic $(n-\nu,\nu+1)$-spacetime} \\
  \Longrightarrow & (M,g) \text{ is a globally hyperbolic $(n-\nu,\nu)$-spacetime} \\
  \Longleftrightarrow &(M', dt^2 \oplus g) \text{ is a globally hyperbolic $(n-\nu+1,\nu)$-spacetime}.
 \end{align*}
 If $\nu = n$ then the first two statements are equivalent.
\end{thm}

 In Example~\ref{ex:notgh} we demonstrate that the first implication in Theorem~\ref{thm:equivsR} \emph{cannot} be reversed if $\nu < n$. Naively speaking, one could even say that this is already evident in the $\nu = 0$ case since \emph{all} manifolds admit a Riemannian metric and \emph{all} Riemannian manifolds are globally hyperbolic (strictly speaking, Riemannian manifolds have empty/degenerate cones and are usually not considered), but already in Theorem~\ref{prop:equivLR} (i) globally hyperbolic Lorentzian products \emph{require} completeness of the Riemannian slice.
 
 In spite of such problems we have shown through our approach that many tools from Lorentzian geometry (cone structures, conformal methods etc.) can successfully be introduced in the semi-Riemannian setting and be utilized to prove interesting new results like Corollary~\ref{semiRcor} generalizing parts of the Hopf--Rinow Theorem. Thus, semi-Riemannian geometry is not only close to Lorentzian and Riemannian geometry in a differential geometric sense (as already emphasized by O'Neill \cite{ON}) but also in a metric context. In fact, it prompts the question what the true Lorentzian (vs.\ general semi-Riemannian) features of recent nonsmooth Lorentzian geometric theories and approaches to quantum gravity are. Answering this question will be crucial for ultimately linking these approaches to the smooth setting of general relativity. Along these lines it could also be insightful to investigate the links between semi-Riemannian and Lorentzian geometry in PDE theory. After all, the notion of global hyperbolicity was introduced in this context by J.\ Leray in 1952 and is crucial for proving global uniqueness of solutions to wave equations.

\medskip
 \textbf{Outline.} In Section~\ref{sec:Lbasics} we recall the basic notions of Lorentzian geometry and causality theory including global hyperbolicity, time functions, and the null distance. We also derive new results for Lorentzian products needed in subsequent sections. In Section~\ref{sec:LRequiv} we show that the Lorentzian  Theorem~\ref{BGHthm} generalizes the metric Hopf--Rinow Theorem~\ref{HRthm} (b)$\Longleftrightarrow$(c) by establishing Theorem~\ref{prop:equivLR} for Lorentzian products. In Section~\ref{sec:cone} we recall basic results of the theory of cone structures and show when the null distance is a well-defined concept in this framework. We also prove the cone completeness-compactness Theorem~\ref{conethm}. In Section~\ref{sec:sR} we introduce the notion of $(n-\nu,\nu)$-spacetime structures for semi-Riemannian manifolds, discuss their existence and links to differential/algebraic topology. We then proceed to show that $(n-\nu,\nu)$-spacetimes are continuous proper cone structures (Theorem~\ref{thm:sRcones}) from which Corollary~\ref{semiRcor} follows. Finally, we relate the notions of (stable) causality and of global hyperbolicity for different $\nu$ and $n-\nu$ for products appearing in Theorem \ref{thm:equivsR}.
 
 \medskip
 \textbf{Notation.} Throughout we denote Riemannian manifolds by $(\Sigma,\sigma)$, Lo\-rentzian and semi-Riemannian manifolds by $(M,g)$. We denote semi-Rie\-mannian products by $(M',g') = (\RR\times M,-dt^2\oplus g)$ and $(M'',g'') = (\RR\times M,dt^2\oplus g)$. If a background Riemannian metric is used it is denoted by $h$. We use the signature convention $(+,\ldots,+)$ for Riemannian manifolds and $(+,\ldots,+,-)$ for Lorentzian manifolds etc.

\section{Causal and metric properties of Lorentzian manifolds}
\label{sec:Lbasics}

 In this section we recall basics notions and results of Lorentzian geometry (for more details see \cites{ON,BEE,Min:LCT}) and the more recent null distance of Sormani and Vega~\cite{SV}. Readers familiar with Lorentzian geometry or the null distance can skip Sections~\ref{ssec:ctime} or \ref{ssec:nulld}, respectively. In Section~\ref{ssec:Lprod} we obtain new results about Lorentzian products that we apply in Section~\ref{sec:LRequiv}.

\subsection{Causality and time}\label{ssec:ctime}

 A \emph{semi-Riemannian manifold} of signature $(n-\nu,\nu)$ is an $n$-dimensional manifold $M$ equipped with metric tensor $g$, i.e., a covariant $2$-tensor field that is symmetric and nondegenerate, with constant index
 \[
  \nu = \max \{ \dim S \, ; \, S \text{ subspace of } T_pM \text{ for which } g_p|_S \text{ is negative definite} \}.
 \]
 A \emph{Lorentzian manifold} has signature $(n-1,1)$. In this section we assume that $(M,g)$ is a Lorentzian manifold, in Section~\ref{sec:sR} we consider the general semi-Riemannian setting.
 
 \medskip
 A tangent vector $v \in T_pM \setminus \{0\}$ is \emph{timelike} if $g(v,v)<0$, \emph{null}  if $g(v,v)=0$, and \emph{spacelike} if $g(v,v)>0$. A tangent vector is \emph{causal} if it is timelike or null. By convention we consider $v=0$ to be spacelike. A Lorentzian manifold $(M,g)$ with a given continuous timelike vector field $X\in\vf(M)$ is said to be time oriented by $X$ and called a \emph{spacetime} \cite{BEE}*{p.\ 25}. It is well known that a manifold admits a spacetime structure if and only if it is either noncompact or compact with vanishing Euler characteristic (see also Introduction).
 
 A causal vector $v \in T_pM \setminus \{0\}$ in a spacetime $(M,g)$ is said to be \emph{future directed} if $g_p(v,X(p)) < 0$. The class of locally Lipschitz curves (with respect to any background Riemannian metric) induces two relations on a spacetime, the \emph{chronological relation}
\begin{align*}
 p \ll q \Longleftrightarrow \exists \text{ future directed timelike curve from } p \text{ to } q,
\end{align*}
 and the \emph{causal relation}
\begin{align*}
 p \leq q \Longleftrightarrow \exists \text{ future directed causal curve from } p \text{ to } q, \text{ or } p=q,
\end{align*}
 with the notation that $p < q$ if $p \leq q$ and $p\neq q$. We write
 \begin{align*}
  I^+(p) &= \{ q\in M \, ; \, p \ll q \}, \\ 
  J^+(p) &= \{ q\in M \, ; \, p \leq q \}
 \end{align*}
 for the chronological future and causal future of $p$, respectively. Analogous definitions apply for the past (with $-$). The chronological relation is open, transitive, and contained in the causal relation, which itself is transitive and reflexive.
 
 Every closed (compact without boundary) Lorentzian manifold admits closed timelike curves, and thus allows time travel. For physical and geometric reasons we would like to exclude this setting and implicitly then focus on the noncompact case.
 
 We require some mild additional property from our spacetimes, namely that they admit time functions.
 
 \begin{defn}\label{def:time}
  Let $(M,g)$ be a spacetime. A continuous function $\tau \colon M \to \RR$ is a \emph{time function} if
  \[
   p < q \Longrightarrow \tau(p) < \tau(q).
  \]
  A smooth time function $\tau$ is called \emph{temporal} if $d\tau(v) > 0$ for all future directed causal vectors $v$ (equivalently, $\nabla\tau$ is past directed timelike).
 \end{defn}

 The existence of a time function conveniently implies that a spacetime is \emph{causal}, i.e., does not admit closed causal curves (and hence the causal relation is antisymmetric, thus an order relation). To be precise, the existence of a time function is equivalent to causality being a stable property, as shown by Hawking~\cite{Haw:time} (see also Proposition~\ref{prop:extime}) who also argues that any physically reasonable theory of gravity must be \emph{stably causal}. The following more refined property of a spacetime is inevitable for a well-posed initial value problem for the Einstein equations in general relativity, the singularity theorems of Penrose and Hawking, Lorentzian splitting theorems, and indeed for most results in Lorentzian geometry.

\begin{defn}
 A spacetime $(M,g)$ is \emph{globally hyperbolic} if it is causal and all causal diamonds $J^+(p) \cap J^-(q)$, $p,q\in M$, are compact.
\end{defn}

 There are many important characterizations and implications of global hyperbolicity which have been obtained since the 1950s, see \cite{BGH2}*{Section 1} for an outline and Theorem~\ref{thm:ghtime} (for proper cone structures) for some characterizations relevant for this work.

\subsection{Null distance}\label{ssec:nulld}

 Sormani and Vega~\cite{SV} showed that any stably causal spacetime $(M,g)$ can be equipped with a conformally invariant length metric space structure that respects and encodes the causal structure and topology of a spacetime. We recall some basic constructions below. See \cites{SV,AB,Ve,SaSo,BGH2,G:null,Al1} for more details and insights, and Section~\ref{ssec:conenull} for some proofs in the more general setting of cone structures.
 
 The chronological relation being open implies that any two points $p,q \in M$ can be joined by a locally Lipschitz path $\beta \colon I \to M$ which is \emph{piecewise causal}, i.e., it is either future or past directed on each subinterval of a partition $\inf I = a_0 < a_1 < \ldots < a_k = \sup I$ of the interval $I$ \cite{SV}*{Lemma 3.5}. The class of piecewise causal curves on $M$ is denoted by $\A$. Since $(M,g)$ is a stably causal spacetime is can be equipped with a time function $\tau$. The \emph{null length} of $\beta \in \A$ is defined as
\[
 \hat L_\tau(\beta) = \sum_{i=1}^k |(\tau \circ \beta)(a_i) - (\tau \circ \beta)(a_{i-1})|.
\]
 It is easy to see that the null length is additive, continuously depends on the curve parameter and is invariant under reparametrizations. Based on this length structure one can obtain a length metric.
 
\begin{defn}
 Let $(M,g)$ be a spacetime with time function $\tau \colon M \to \RR$. The \emph{null distance} $\hat d_\tau$ between to points $p,q \in M$ is
 \[
  \hat d_\tau(p,q) = \inf \{ \hat L_\tau(\beta) \, ; \, \beta \in \A \text{ from } p \text{ to } q \}.
 \]
\end{defn}

 By definition, $\hat d_\tau$ is a conformally invariant pseudometric on $(M,g)$. That the null distance is distinguishing asks slightly more of the time function used. Such time functions always exist already in the stably causal setting, for example, one may use temporal functions~\cite{BeSa:orthsplit}.
 
\begin{thm}[Sormani--Vega~\cite{SV}*{Theorem 4.6}]\label{thm:defnulldist}
 Let $(M,g)$ be a spacetime with \emph{locally anti-Lipschitz} time function $\tau \colon M \to \RR$, i.e., for every point there exist a neighborhood $U$ and Riemannian metric $h$ such that
 \[
  p,q \in U \text{ and } p \leq q \Longrightarrow \tau(q) - \tau(p) \geq d_h(p,q).
 \]
 Then the null distance $\hat d_\tau$ is a conformally invariant metric on $M$ that induces the manifold topology.
\end{thm}

 The author and B.\ Allen have shown the following important connection to completeness that was used in the proof of Theorem~\ref{BGHthm} (c')$\Longrightarrow$(b').
 
 \begin{thm}[Allen--Burtscher~\cite{AB}*{Theorem 1.6, Corollary 3.15}]\label{ABthm}
  Let $(M,g)$ be a spacetime with time function $\tau$ that is globally anti-Lipschitz with respect to a complete (Riemannian) metric. Then $(M,\hat d_\tau)$ is a complete metric space.
 \end{thm}
 
 Note that the completeness needed in the assumption of the statement is only the metric completeness property (b) in the Hopf--Rinow Theorem~\ref{HRthm}, in particular, the proof does \emph{not} require the use of the infinity postulate (a), as can be seen from the fact that it holds for \emph{any} complete metric that induces the manifold topology~\cite{AB}*{Theorem 1.6}. The author and L.\ Garc\'ia-Heveling have in \cite{BGH2}*{Section 1} established more connections of the null distance to Riemannian geometry. It should be noted, in particular, that weak temporal functions are precisely those time functions $\tau$ that induce locally equivalent metrics $\hat d_\tau$ \cite{BGH2}*{Corollary 1.8}.

\subsection{Lorentzian products}\label{ssec:Lprod}

 Every Riemannian manifold $(\Sigma,\sigma)$ can be studied as a Lorentzian manifold by considering the canonical product
 \[
  (\Sigma',\sigma') = (\RR \times \Sigma, -dt^2 \oplus \sigma).
 \]
 In what follows we prove some results about such Lorentzian product spacetimes needed in Section~\ref{sec:LRequiv}. Although our focus is on the above case it is easy to see that most results also hold for products $N \oplus \Sigma$ for $N$ a suitable Lorentzian manifold and that similar results are true also for warped product spacetimes (see also \cites{AB,BEE}).
 
 \medskip
 The canonical projection $t = \pi_0 \colon \Sigma' \to \RR$, $t(p_0,p_\Sigma) := p_0$, satisfies
\[
 dt_p(v) = - g(v,\partial_t) = v_0 > 0
\]
 for future directed causal vectors $v \in T_p\Sigma' \cong \RR \times T_{p_\Sigma}\Sigma$. Thus $t$ is a smooth temporal function of $(M,g)$. The null distance $\hat d_t$ on $\Sigma'$ is thus well-defined by Theorem~\ref{thm:defnulldist} and by \cite{AB}*{Lemma 4.4} of the form 
\begin{align}\label{dtau}
 \hat d_t(p,q) = \begin{cases}
                     |t(q) - t(p)|, & q \in J^\pm(p), \\
                     d_\sigma(p_\Sigma,q_\Sigma), & \text{otherwise}.
                    \end{cases}
\end{align}
We show slightly more in the following result.

\begin{lem}\label{lem1}
 Let $(\Sigma,\sigma)$ be a Riemannian manifold and $(\Sigma',\sigma')$ be the canonical Lorentzian product. Then
 \begin{align}\label{dtmax}
  \hat d_t(p,q) = \max \{ |t(q)-t(p)|, d_\sigma(p_\Sigma,q_\Sigma) \},
 \end{align}
 where $t(p) = \pi_0(p) = p_0$ and $p_\Sigma = \pi_\Sigma(p)$ are the orthogonal projections onto $\RR$ and $\Sigma$, respectively.
\end{lem}

\begin{proof}
 By \eqref{dtau} it is clear that the $\leq$ inequality in \eqref{dtmax} holds. It remains to be shown that the converse $\geq$ inequality holds too. We distinguish the two cases mentioned in \eqref{dtau}.
 
 If $q \in J^+ (p)$, then there exists a future directed causal curve $\gamma \colon I \to M$ from $p$ to $q$. We write the different components of $\gamma$ as $\gamma_0 = t\circ \gamma$ and $\gamma_\Sigma = \pi_\Sigma \circ \gamma$. That $\gamma$ is causal means that $g(\dot\gamma,\dot\gamma) \leq 0$, hence $\dot \gamma_0 = |dt(\dot\gamma)| \geq \| \dot\gamma_\Sigma\|_g$ (the positive sign of $\dot\gamma_0$ is due to the future directedness of $\gamma$),
 and thus by the fundamental theorem of calculus
 \[
  t(q)-t(p) = \int_I \dot\gamma_0(s) ds \geq \int_I \|\dot\gamma_\Sigma(s)\|_g ds = L_g(\gamma_\Sigma) \geq d_\sigma(p_\Sigma,q_\Sigma).
 \]
 The proof for $q \in J^-(p)$ is the same (just exchange the role of $p$ and $q$). Thus in case $p$ and $q$ are causally related the reverse inequality $\geq$ in \eqref{dtmax} holds. 

 Suppose that $q \not \in J^\pm(p)$. Then it follows immediately from the definition of the null distance and \eqref{dtau} that
\[
 d_\sigma(p_\Sigma,q_\Sigma) = \hat d_t(p,q) \geq |t(q)-t(q)|,
\]
 and we have thus shown $\geq$ of \eqref{dtmax}.
\end{proof}

 We use Lemma~\ref{lem1} to characterize the boundary of the causal future.

\begin{lem}\label{lem2}
 Let $(\Sigma,\sigma)$ be a Riemannian manifold and let $(\Sigma',\sigma')$ be the canonical Lorentzian product. Then
  \begin{align}\label{eq:bd}
  q \in \partial J^+(p) \Longleftrightarrow t(q)-t(p) = \hat d_t(p,q) = d_\sigma(p_\Sigma,q_\Sigma).
 \end{align}
\end{lem}

\begin{proof}
 ($\Longrightarrow$) Suppose that $q \in \partial J^+(p) = \overline{J^+(p)} \cap \overline{J^+(p)^\complement}$. The identity \eqref{dtau} together with the continuity of $t$ and $\hat d_t$ implies the first equality in \eqref{eq:bd} for all $q \in \overline{J^+(p)}$ (approximation from inside), and \eqref{dtau} together with the continuity of $\pi_\Sigma$, $d_\sigma$ and $\hat d_t$ implies the second equality for all $q \in \overline{J^+(p)^\complement}$ (approximation from outside).
 
 ($\Longleftarrow$) That $q$ is a boundary point if both equalities in \eqref{eq:bd} hold follows from Lemma~\ref{lem1} in a similar fashion.  To be more precise, we show that
 \[
  d_\sigma(p_\Sigma,q_\Sigma) = t(q) - t(p)
 \]
 leads to a contradiction if $q$ is not a a boundary point of $J^+(p)$ based on openness and the product structure:
 
 Suppose $q \in I^+(p)$. Since $I^+(p)$ is open and $\hat d_t$ induces the manifold topology, there exists an $\varepsilon >0$ such that
 \[
  \hat B_{2\varepsilon}(q) = \{ x \in \Sigma' \, ; \, \hat d_t(q,x) < 2 \varepsilon \} \subseteq I^+(p).
 \]
Clearly, the point $q' = (q_0-\varepsilon,q_\Sigma) \in \hat B_{2\varepsilon}(q)$ and thus satisfies
 \[
  t(q') - t(p) = t(q)-\varepsilon - t(p) < d_\sigma(p_\Sigma,q_\Sigma) = d_\sigma(p_\Sigma,q'_\Sigma).
 \]
 By \eqref{dtau} the left hand side is equal to $\hat d_t(p,q')$, hence this contradicts \eqref{dtmax} applied to $p$ and $q'$.
 
 Suppose $q \in \overline{I^+(p)}^\complement$. Since $\overline{I^+(p)}^\complement$ is open, similarly, we have a point $q'' = (q_0 + \varepsilon,q_\Sigma) \in \overline{I^+(p)}^\complement$ satisfying
 \[
  t(q'') - t(p) = t(q)+\varepsilon - t(p) > d_\sigma(p_\Sigma,q_\Sigma) = d_\sigma(p_\Sigma,q'_\Sigma).
 \]
 By \eqref{dtau} the right hand side is equal to $\hat d_t(p,q'')$, again a contradiction to \eqref{dtmax}.
 
 Because of the disjoint union $\Sigma' = I^+(p) \sqcup \partial I^+(p) \sqcup (\overline{I^+(p)})^\complement$, and the fact that we already know that the equalities in \eqref{eq:bd} hold \emph{at least} on the set $\partial I^+(p) = \partial J^+(p)$ by the first implication ($\Longrightarrow$), we are done.
\end{proof}

 The above findings are interesting independent of completeness or the Hopf-Rinow Theorem. For us they are crucial for proving Theorem~\ref{prop:equivLR} in Section~\ref{sec:LRequiv} because we do \emph{not} assume geodesic completeness of $(\Sigma,\sigma)$ nor any knowledge of the Hopf--Rinow Theorem~\ref{HRthm} in parts or as a whole (otherwise we would obtain a circular argument). It is useful to note, however, that with the additional assumption of geodesic completeness a lot more can be shown. 

\begin{rem}[Geodesic completeness of the fiber and global hyperbolicity]\label{rem:warped}
 It is long known that a Lorentzian warped product $I \times_f \Sigma$ is globally hyperbolic if and only if the fibers $(\Sigma,\sigma)$ are geodesically complete Riemannian manifolds \cite{BEE}*{Theorem 3.66}. In the case of Lorentzian products one can show that both statements are also equivalent to the Lorentzian product $(I \times \Sigma, -dt^2 \oplus \sigma)$ being geodesically complete \cite{BEE}*{Theorem 3.67}.
\end{rem}

 Geodesic completeness of $(\Sigma,\sigma)$ can furthermore be employed constructively for proofs about the null distance on $\RR \times \Sigma$.

\begin{rem}[Geodesic completeness of the fiber and causality encoding of the null distance]
 Suppose $\Sigma$ is a geodesically complete Riemannian manifold. Given $p,q \in \Sigma'=I \times\Sigma$ we can then always construct a length-minimizing Riemannian geodesic $\alpha$ between $p_\Sigma$ and $q_\Sigma$ in $\Sigma$. We can lift $\alpha$ to a curve $\gamma(s) = (s,\alpha(s))$ (recall that geodesics satisfy $\|\dot\alpha\|_\sigma = \text{const.}$) connecting $\gamma(t(p))=p$ to $\gamma(t(q)) = q$ in $M$. This is particularly useful when $p$ and $q$ are such that
 \[
  \hat d_t(p,q) = t(q) - t(p)
 \]
 holds, because by Lemma~\ref{lem1} then $L_\sigma(\alpha) = d_\sigma(p_\Sigma,q_\Sigma) \leq \hat d_t(p,q) = t(q)-t(p)$ and thus $\|\dot\alpha\|_\sigma \leq 1$. In other words, $\| \dot \gamma\|_g = -1 + \|\dot\alpha\|_\sigma \leq 0$, meaning that $\gamma$ is a future directed causal curve from $p$ to $q$. Hence we have shown that the canonical null distance of the Lorentzian warped product $(I \times \Sigma, -dt^2 \oplus \sigma)$ encodes causality, i.e.,
 \begin{align}\label{encc}
  p \leq q \Longleftrightarrow \hat d_t(p,q) = t(q) - t(p).
 \end{align}
 Of course we could have also inferred this from Lemma~\ref{lem2} and the definition of the null distance (since ($\Longrightarrow$) in \eqref{encc} holds always). With a significantly more involved proof one can show that Lorentzian warped products with complete fibers as well as all globally hyperbolic spacetimes encode causality globally for sensible time functions (see \cite{SV}*{Theorem 3.25} and \cite{BGH2}*{Theorem 1.9}).
\end{rem}

 On the other hand, it is easy to construct examples of Lorentzian products with incomplete fibers that do not satisfy \eqref{encc}. Consider, for example, $\RR \times (\RR^n \setminus \{0\})$ equipped with the metric induced from Minkowski metric and the canonical time function. Nevertheless, any null distance on any spacetime encodes causality locally as has been shown independently by Sakovich and Sormani~\cite{SaSo}*{Theorem 1.1} and Garc\'ia-Heveling and the author~\cite{BGH2}*{Theorem 3.4}.

\begin{thm}[\cite{BGH2}*{Theorem 3.4}]\label{thm:localcausality}
 Let $(M,g)$ be a spacetime and $\tau$ a temporal function. Then at every point $x \in M$ there is a neighborhood $U$ of $x$ such that for all $p,q\in U$
 \[
  p \leq q \Longleftrightarrow \hat d_\tau(p,q) = \tau(q) -\tau(p).
 \]
\end{thm}

\section{Theorem~\ref{BGHthm} implies the metric Hopf--Rinow Theorem}\label{sec:LRequiv}

 In this section we are solely concerned with Lorentzian products, i.e., spacetimes of the form \[(\Sigma',\sigma')=(\RR \times\Sigma,-dt^2\oplus \sigma)\]
 with $(\Sigma,\sigma)$ a Riemannian metric. We prove Theorem~\ref{prop:equivLR} (i) and (ii) in Sections~\ref{ssec:HBGH} and \ref{ssec:complete} respectively. In other words we show that the metric Hopf--Rinow Theorem~\ref{HRthm} (b)$\Longleftrightarrow$(c) is a special case of Theorem~\ref{BGHthm}, and thus can be also viewed as a special case of Theorem~\ref{conethm}.

\subsection{The Heine--Borel property vs.\ global hyperbolicity}\label{ssec:HBGH}
 
 The product $(\Sigma',\sigma')$ is always a \emph{stably causal} spacetime independent of properties of $(\Sigma,\sigma)$, thus it remains to relate the Riemannian compactness postulate to the corresponding conformal compactness postulate for Lorentzian products.

 Strictly speaking we prove nothing new, but our proof is new: It is well-known that geodesic completeness of the Riemannian fiber $(\Sigma,\sigma)$ is equivalent to be the Lorentzian product $I \times \Sigma$ being globally hyperbolic (see Remark~\ref{rem:warped}), and together with the Hopf--Rinow Theorem~\ref{HRthm} (a)$\Longleftrightarrow$(b) we would therefore be done. But we need to abolish the sophisticated analytic notion of geodesic completeness and just relate the two topological conditions directly. By doing away with (geodesic) completeness we can avoid a circular argument and thus ensure that Theorem~\ref{BGHthm} really implies Theorem~\ref{HRthm} (b)$\Longleftrightarrow$(c).

\begin{prop}\label{prop1}
 Let $(\Sigma,\sigma)$ be a Riemannian manifold and let $(\Sigma',\sigma')=(\RR \times \Sigma,-dt^2 \oplus \sigma)$ be the corresponding Lorentzian product manifold. Then
 \[ (\Sigma,\sigma) \text{ has the Heine--Borel property } \Longleftrightarrow (\Sigma',\sigma') \text{ is globally hyperbolic}.
 \]
\end{prop}

\begin{proof}
 By definition, $t$ is a temporal function for $\Sigma'$, thus $\Sigma'$ is stably causal. it thus remains to be shown that
 \begin{align*}
  (\Sigma,\sigma) \text{ has the } &\text{Heine--Borel property } \\ &\Longleftrightarrow \text{ causal diamonds in } (\Sigma',\sigma') \text{ are compact}.
 \end{align*}

 $(\Longrightarrow)$ Suppose $(\Sigma,\sigma)$ satisfies the Heine--Borel property. Let $p,q \in \Sigma'$ and consider the causal diamond $J^+(p) \cap J^-(q)$. If $J^+(p) \cap J^-(q) = \emptyset$ nothing is to prove, so suppose that $J^+(p) \cap J^-(q) \neq \emptyset$. Due to the definiteness of $\hat d_t$ it follows that $\hat d_t(p,q) = t(q) - t(p) >0$ due to the definiteness of $\hat d_t$. We show that $J^+(p) \cap J^-(q)$ is (i) contained in a compact set and (ii) closed, hence compact itself. 
 
 (i) Let $r := t(q)-t(p)>0$. For any $x \in J^+(p) \cap J^-(q)$, by \eqref{dtau},
 \begin{align}\label{test}
  t(x) \in [t(p),t(q)].
 \end{align}
 Moreover, \eqref{test} together with Lemma~\ref{lem1} implies that
 \begin{align*}
  \max\{ d_\sigma(p_\Sigma,x_\Sigma), d_\sigma(x_\Sigma,q_\Sigma) \} &\leq \max\{ d_t(p,x), d_t(x,q) \} \leq r,
 \end{align*}
 and hence
 \begin{align}\label{gest}
  x_\Sigma \in \overline{B_r^\sigma(p_\Sigma) \cap B_r^\sigma(q_\Sigma)}.
 \end{align}
 Together, \eqref{test} and \eqref{gest} imply that the causal diamond is contained in the intersection of two closed cylinders,
 \[
  J^+(p) \cap J^-(q) \subseteq [t(p),t(q)] \times \overline{B_{r}^\sigma(p_\Sigma) \cap B_{r}^\sigma(q_\Sigma)} =: A(p,q).
 \]
 The interval $[t(p),t(q)]$ is closed and hence compact in $\RR$. Moreover, the set $\overline{B_r^\sigma(p_\Sigma) \cap B_{r}^\sigma(q_\Sigma)}$ is bounded and closed in $\Sigma$, thus by the Heine--Borel property of $\Sigma$ it is also compact. Hence $A(p,q)$ is compact.
 
 (ii) It remains to be shown that the causal diamonds $J^+(p) \cap J^-(q)$ are closed. We show that $J^+(p)$ is closed or, to be more precise, that $\partial J^+(p) \subseteq J^+(p)$. Suppose that $x = (x_0,x_\Sigma) \in \partial J^+(p)$ with $x_0 = t(x)$. We construct a particular sequence $(x^n)_n$ approximating $x$ from inside. By Lemma~\ref{lem2} 
 \[
  T:= \hat d_t(p,x) = t(x) - t(p) = d_\sigma(p_\Sigma,x_\Sigma).
 \]
 In particular, $x_\Sigma \in \partial B^\sigma_T(p_\Sigma) \subseteq \Sigma$, meaning that there is a sequence $x_\Sigma^n \in B^\sigma_T(p_\Sigma) \subseteq \Sigma$ in the interior of this ball approximating $x_\Sigma$, i.e., satisfying
 \[
  d_\sigma(x^n_\Sigma,x_\Sigma) < \frac{1}{n}, \qquad d_\sigma(x_\Sigma^n,p_\Sigma) < T,
 \]
 We can lift the sequence $(x^n_\Sigma)_n$ to a sequence $(x^n)_n$ in $\Sigma'$ by setting $x^n = (x_0,x^n_\Sigma)$. Then $x^n$ converges to $x$ with respect to the null distance $\hat d_t$ because
 \[
  \hat d_t(x,x^n) = \max \{0, d_\sigma(x^n_\Sigma,x_\Sigma)\} < \frac{1}{n}
 \]
 by Lemma \ref{lem1}. Since $d_\sigma(x_\Sigma^n,p_\Sigma) < T = \hat d_t(x^n,p)$ it follows from \eqref{dtau} and Lemmas~\ref{lem1} and \ref{lem2} that $x^n \in I^+(p) = J^+(p) \setminus \partial J^+(p)$, hence there exists a timelike curve $\gamma_n$ with null length
 \[
 \hat L_t(\gamma_n) = t(x^n) - t(p) = t(x) - t(p) = T.
 \]
 Moreover, each $\gamma_n$ is contained in the compact set $A(p,q)$ of Step (i). Recall that Allen and the author have shown that for piecewise causal curves $\hat L_t = L_{\hat d_t}$ (the latter being lower semicontinuous by \cite{BBI}*{Prop.\ 2.3.4} for rectifiable paths) and that $\hat d_t$ is an intrinsic metric \cite{AB}*{Prop.\ 3.8}. By the Arzela--Ascol\'i Theorem a subsequence of $(\gamma_n)_n$ thus uniformly converges to a rectifiable limit curve $\gamma$ in $A(p,q)$ connecting $p$ and $x$, which is Lipschitz with respect to $\hat d_t$ (since the image of $\gamma$ is compact, it is even Lipschitz continuous in the usual sense with respect to any Riemannian metric by \cite{BGH2}*{Theorem 1.7}) and by lower semicontinuity of the length functional
 \begin{align}\label{Lhatdt}
  t(x) - t(p) = \hat d_t(p,x) \leq L_{\hat d_t}(\gamma) \leq \lim_{n\to\infty} \hat L_t(\gamma_n) = t(x) - t(p) = T.
 \end{align}
 Assume without loss of generality that $\gamma \colon [0,1] \to M$. Condition \eqref{Lhatdt} can be localized, meaning that
 \begin{align}\label{Lhatlocal}
  L_{\hat d_t}(\gamma|_{[s,s']}) \leq t(\gamma(s')) - t(\gamma(s)),
 \end{align}
 because \eqref{Lhatdt} and the converse of \eqref{Lhatlocal} would imply
 \begin{align*}
  t(x) - t(p) &= \hat d_t(p,x) = L_{\hat d_t}(\gamma) \\
  &> |t(x) - t(\gamma(s'))| + t(\gamma(s'))-t(\gamma(s)) + |t(\gamma(s)) - t(p)| \\ &\geq t(x) - t(p),
 \end{align*}
 a contradiction.
 
 It remains to be shown that $\gamma$ connects $p$ and $x$ causally. Consider the set
 \[
  A := \{ s \in [0,1] \, ; \, \gamma(s) \in J^+(p) \}.
 \]
 Clearly, $A \neq \emptyset$ because $\gamma(0)= p \in J^+(p)$. By \eqref{Lhatlocal} all $s,s' \in [0,1]$ satisfy
\begin{align}\label{dtequal}
 \hat d_t(\gamma(s),\gamma(s')) = t(\gamma(s')) - t(\gamma(s)).
\end{align}
 Since the null distance locally around $p$ encodes causality by \cite{BGH2}*{Theorem 3.4}, we thus know that $\gamma(s) \in J^+(p)$ for small $s$, i.e., $A$ is open. Suppose
 \[
  s_0 := \sup A.
 \]
 Again, local causality encoding around $\gamma(s_0)$ together with \eqref{dtequal} applied to $s \in A$ sufficiently close to $s_0 > s$, yields that $s_0 \in A$, i.e., $A$ is also closed. Thus $A=[0,1]$ and therefore $x =\gamma(1) \in J^+(p)$, as desired. 
 
 \medskip
 $(\Longleftarrow)$ Suppose $C$ is a closed and bounded subset of $\Sigma$, and the causal diamonds in $(\Sigma',\sigma')$ are compact. The product $\widetilde C := \{0\} \times C$ is closed in $\Sigma'$. We can assume that $C$ is nonempty, and that there is a $p_\Sigma \in \Sigma$ and $r>0$ such that $C$ is contained in the open $d_\sigma$-ball $B_r^\sigma(p_\Sigma)$. 
 
 If $x \in C$, then $d_\sigma(p_\Sigma,x) < r$, and hence there is a curve $\alpha$ from $p_\Sigma$ to $x$ in $\Sigma$ such that $L_\sigma(\alpha) < r$. We may assume that $\alpha$ is parameterized by constant speed $\|\dot\alpha\|_\sigma \leq 1$ on $[-r,0]$, and lift it to a curve $\gamma(s):= (s,\alpha(s))$ from $p^- := (-r,p_\Sigma)$ to $(0,x)$ in $\Sigma'$. Since
 \[
  \sigma'(\dot\gamma,\dot\gamma) = -1 + \| \dot\alpha\|_\sigma \leq 0,
 \]
 $\gamma$ is future directed causal and hence $(0,x) \in J^+(p^-)$. Similarly it follows that $(0,x) \in J^-(p^+)$ for $p^+ := (r,p_\Sigma)$. Thus the closed set $\widetilde C$ is contained in the (by assumption) compact causal diamond $J^+(p^-) \cap J^-(p^+)$, and therefore compact itself. Since the canonical projection $\pi_\Sigma \colon \Sigma' \to \Sigma$ onto the second component is continuous (in the product topology), $C$ is also be compact. Therefore, $\Sigma$ satisfies the Heine--Borel property.
\end{proof}

\subsection{Metric completeness}\label{ssec:complete}

 In the same spirit as in Section~\ref{ssec:HBGH} we directly relate the metric completeness property of $\Sigma$ to that of $\Sigma' = \RR \times \Sigma$.

\begin{prop}\label{prop2}
 Let $(\Sigma,\sigma)$ be a Riemannian manifold and $(\Sigma',\sigma') = (\RR \times \Sigma, -dt^2 \oplus \sigma)$ be the canonical Lorentzian product manifold with time function $t = \pi_0$. Then the metric spaces satisfy
 \[
  (\Sigma,d_\sigma) \text{ is complete } \Longleftrightarrow (\Sigma',\hat d_t) \text{ is complete}.
 \]
\end{prop}

\begin{proof}
 $(\Longrightarrow)$ Assume that $(\Sigma,d_\sigma)$ is a complete metric space. To show the same of $(\Sigma',\hat d_t)$ assume that $p_n = (p^n_0,p^n_\Sigma)_n$ is a Cauchy sequence in $(\Sigma',\hat d_t)$. Since by Lemma~\ref{lem1}
 \[
  \hat d_t (p_m,p_n) \geq \frac{1}{2} \left(|p^m_0-p^n_0| + d_\sigma(p^m_\Sigma,p^n_\Sigma) \right)
 \]
 it follows that $(p^n_0)_n$ and $(p^n_\Sigma)_n$ are Cauchy sequences in $\RR$ and $\Sigma$, respectively. Since both spaces are complete, there is an $p_\infty = (p^\infty_0,p^\infty_\Sigma) \in \RR \times \Sigma$ such that, again by Lemma~\ref{lem1}, as $n \to\infty$,
 \[
  \hat d_t (p_n,p_\infty) = \max \{ |p^n_0-p^\infty_0|,d_\sigma(p^n_\Sigma,p^\infty_\Sigma) \} \to 0.
 \]
 Hence $p_\infty$ is the limit of $(p_n)_n$, and $(\Sigma',\hat d_t)$ is complete.
 
 \medskip
 $(\Longleftarrow)$ Suppose $(\Sigma',\hat d_t)$ is a complete metric space. Let $(p^n_\Sigma)_n$ be a Cauchy sequence in $\Sigma$. Since by \eqref{dtau}
 \[
  \hat d_t((0,p^m_\Sigma),(0,p^n_\Sigma)) = d_\sigma(p^m_\Sigma,p^n_\Sigma)
 \]
 the sequence $p_n = (0,p^n_\Sigma)_n$ is a Cauchy sequence in $\Sigma'$. Hence there exists a limit point $p_\infty = (p^\infty_0,p^\infty_\Sigma)$ in $\Sigma'$. 
 Since $t = \pi_0$ is continuous it follows that $p^\infty_0=0$. Therefore, as $n \to \infty$,
 \[
  d_\sigma(p^n_\Sigma,p^\infty_\Sigma) = \hat d_t(p_n,p_\infty) \to 0,
 \]
 and so $p^\infty_\Sigma$ is the limit point of the sequence $(p^n_\Sigma)_n$ in $(\Sigma,d_\sigma)$. Thus $(\Sigma,d_\sigma)$ is complete.
\end{proof}

\section{Cone structures}
\label{sec:cone}

 The aim of this section is to introduce all notions of and to prove Theorem~\ref{conethm}. In Sections~\ref{ssec:cones} and \ref{ssec:conetimes} we recall some definitions and results of the theory of cone structures. Readers familiar with the recent theory of cone structures as in \cites{Min:cone,FaSi,Fa:time,BeSu:time,BeSu:gh} may skip this part. In Section~\ref{ssec:conenull} we introduce the null distance and show that it is a well-defined metric on proper cone structures. In Section~\ref{ssec:coneHR} we prove that the Lorentzian Theorem~\ref{BGHthm} also has a cone version, namely Theorem~\ref{conethm}. These extensions are straightforward but nonetheless have to be checked carefully. The main asset of setting up everything in the cone framework is, apart from its own benefit, that we can efficiently apply it to our new notion of semi-Riemannian spacetimes in Section~\ref{sec:sR}.

\subsection{Cone structures}\label{ssec:cones}

 We use mostly the conventions of Minguzzi~\cite{Min:cone} which builds on the pioneering work of Fathi and Siconolfi~\cites{FaSi,Fa:time} and results of Bernard and Suhr~\cites{BeSu:time,BeSu:gh}.

\begin{defn}\label{def:vectorcone}
 Let $V$ be a finite-dimensional vector space. A \emph{convex cone} $C$ is a convex subset of $V \setminus \{0\}$ such that
 \[
  v \in C, s >0 \Longrightarrow sv \in C.
 \]
 A cone $C$ is \emph{sharp} (or regular \cite{BeSu:time}) if $C \cup \{0\}$ does not contain any line passing through the origin. A cone is called \emph{closed} if it is closed in the topology induced on $V \setminus \{0\}$ by the topology of $V$.
 
 A \emph{proper cone} is a closed sharp convex cone with nonempty interior.
\end{defn}

\begin{defn}\label{def:constructure}
 Let $M$ be a manifold. A \emph{cone structure} $(M,C)$ is a multivalued map $p \mapsto C_p$, where $C_p \subseteq T_pM \setminus \{0\}$ is a closed sharp convex nonempty cone.
 
 A \emph{closed cone structure} is a cone structure which is a closed subbundle of the slit tangent bundle $TM \setminus (TM)_0$, where $(TM)_0 = \{ 0_p \, ; \, 0 \in T_pM\}$ denotes the zero section.

 A \emph{proper cone structure} is a closed cone structure in which the cone bundle is proper, i.e., $(\inte C)_p \neq \emptyset$ for all $p \in M$.
\end{defn}

 Note that for a closed or proper cone structure it is \emph{not} sufficient to interpret these notions fiberwise, i.e, that each cone $C_p$ is closed or proper! Instead, these notions require the use of the topology of the cone bundle. Higher regularity generally helps. For many properties, approximation by wider and more regular cones is sufficient. We write $C \prec C'$ if $C \subseteq \inte C'$ and $C \preccurlyeq C'$ if $C \subseteq C'$.

\begin{prop}[\cite{Min:cone}*{Proposition 2.6}]
 For a $C^0$ proper cone structure $(\inte C)_p = \inte C_p$ for every $p$.
\end{prop}

\subsection{Causality and time functions}\label{ssec:conetimes}

By using the differential inclusions
 \[
  \dot\gamma(t) \in C_{\gamma(t)}, \qquad \dot\gamma(t) \in (\inte C)_{\gamma(t)}
 \]
 on a closed or proper cone structure $(M,C)$ one can define a chronological relation $I^+$ and a causal relation $J^+$ on $(M,C)$ that resembles the corresponding Lorentzian relations (cf.\ Section~\ref{ssec:ctime}). In addition to recalling the notation of~\cite{Min:cone}*{Section 2.1} we add some elements. For various reasons it is convenient and sufficient to work within the class of locally absolutely continuous (or locally Lipschitz) paths \cites{Bur,Min:cone}.
 
\begin{defn}\label{def:causalcurve}
 Let $(M,C)$ be a closed cone structure. An element $v \in C_p \subseteq T_pM$ is called \emph{future directed causal vector}. An element $w \in (\inte C)_p \subseteq T_pM$ is called \emph{future directed timelike vector}. A vector $z \in (\partial C)_p = C_p \setminus (\inte C)_p$ is called \emph{future directed lightlike}. A vector $v$ is called \emph{past directed} causal/timelike/lightlike if $-v$ is future directed causal/timelike/lightlike.
 
 A \emph{future directed causal curve} is the image of an absolutely continuous solution to the differential inclusion $\dot\gamma(t) \in C_{\gamma(t)}$. A \emph{future directed timelike curve} is the image of a (piecewise) $C^1$ solution of $\dot\gamma(t) \in (\inte C)_{\gamma(t)}$. The definitions for past directed curves are analogous.
\end{defn}
 
 Based on Definition~\ref{def:causalcurve} one proceeds as in Section~\ref{ssec:ctime} to define the \emph{causal relation} $J^+$ and \emph{chronological relation} $I^+$. The usual properties follow. We review these results briefly to remind the reader that (and why) several require proper cone structures as opposed to only closed ones.

\begin{prop}[\cite{Min:cone}*{Proposition 2.8}]\label{prop:IJ}
 Let $(M,C)$ be a closed cone structure. The causal relation $J^\pm$ is transitive and reflexive, i.e., a preorder. The chronological relation $I^\pm$ is transitive and contained in $J^\pm$. If $(M,C)$ is proper then $I^\pm$ is open and nonempty.
\end{prop}

\begin{proof}[Sketch of proof]
 Clearly, $I^\pm \subseteq J^\pm$.
 Reflexivity of $J^\pm$ is clear by definition. Transitivity is clear for both relations since it just requires concatenating paths. That $I^\pm$ is open follows from \cite{Min:cone}*{Proposition 2.8} since the corresponding cone structure $C$ is proper and one can then use an approximation by a $C^0$ proper cone structure. Furthermore, that $I^\pm$ is nonempty follows from the assumption that there exists $v \in (\inte C)_p \neq \emptyset$ for all $p$ and the fact that there is a timelike curve through $p$ with velocity $v$ \cite{Min:cone}*{Theorem 2.2}.
\end{proof}

\begin{cor}\label{cor:closedcurves}
 A closed manifold with proper cone structure $(M,C)$ admits closed timelike curves.
\end{cor}

\begin{proof}
 The collection $\{I^+(p) \, ; \, p \in M\}$ is an open cover of $M$ by Proposition \ref{prop:IJ}. Since $M$ is compact there exists a finite subcover $\{I^+(p_i) \, ; \, i =1,\ldots,m\}$. Suppose $p_i \not\in I^+(p_i)$ for all $i$. Then $p_1 \in I^+(p_{i_1})$ for some $p_{i_1} \neq p_1$ and there is a future directed timelike curve from $p_{i_1}$ to $p_1$, and subsequently from $p_{i_2}$ to $p_{i_1}$, $i_1 \neq i_2$, etc. Since only finitely many $p_i$ exist after at most $m+1$ steps a repetition occurs, and hence there exists a closed timelike curve.
\end{proof}

 In order to extend the notion of null distance to semi-Riemannian manifolds, we recall that time and temporal functions for closed cone structures are defined in the same way as in Definition~\ref{def:time}. One can also define locally anti-Lipschitz time functions in the same way as in Theorem~\ref{thm:defnulldist}. See also \cite{Min:cone}*{Section 2.2} for more notions.

 \medskip
 As in Lorentzian geometry the existence of time functions trivially requires that $(M,C)$ is causal. Fathi and Siconolfi \cite{FaSi}*{Theorem 1.1} were the first to show that stable causality is sufficient. Minguzzi's proof is based on the $K^+$-relation of Sorkin and Woolgar \cite{SoWo} and Levin's Theorem~\cite{Lev}.

\begin{prop}[\cite{Min:cone}*{Theorem 2.30}]\label{prop:extime}
 Let $(M,C)$ be a closed cone structure. The following are equivalent:
 \begin{enumerate}
  \item $(M,C)$ is \emph{$K$-causal}, i.e., the smallest closed and transitive relation $K^+$ containing $J^+$ is antisymmetric,
  \item $(M,C)$ is \emph{stably causal}, i.e., there is a $C^0$ proper cone structure $C'$ such that $C \prec C'$ which is causal,
  \item there exists a time function on $M$,
  \item there exists a smooth temporal function on $M$.
 \end{enumerate}
\end{prop}

%

 Several important characterizations of global hyperbolicity are already known for closed and proper cone structures as well as the stability of this notion. For continuous cone structures these results were obtained in the pioneering work of Fathi and Siconolfi~\cite{FaSi}. We only state the results that we will use and refer to the literature \cites{FaSi,Fa:time,BeSu:time,BeSu:gh,Min:cone} for proofs.

\begin{defn}\label{def:ggh}
 A closed cone structure $(M,C)$ is called \emph{globally hyperbolic} if it is causal and all causal diamonds $J^+(p) \cap J^-(q)$, $p,q\in M$, are compact.
\end{defn}

\begin{thm}[\cite{Min:cone}*{Theorem 2.45}]\label{thm:ghtime}
 Let $(M,C)$ be a closed cone structure. Then the following properties are equivalent:
 \begin{enumerate}
  \item $(M,C)$ is globally hyperbolic,
  \item there exists a \emph{Cauchy time function}, i.e., a time function $\tau \colon M \to \RR$ such that for every inextendible future/past directed causal curve $\gamma$ we have $(\tau \circ \gamma)(\RR) = \RR$,
  \item there exists a smooth \emph{completely uniform temporal function} (automatically Cauchy) on $M$, i.e., a temporal function $\tau \colon M \to \RR$ for which there exists a complete Riemannian metric $h$ such that $d\tau(v) \geq \| v\|_h$ for all future directed causal vectors $v$,
  \item there exists a (stable) \emph{Cauchy hypersurface}, i.e., an acausal (no two points are connected by a future/past directed causal curve) topological hypersurface $\Sigma$ such that $D(\Sigma) = D^+(\Sigma)\cup D^-(\Sigma) = M$, where
   \begin{align*}
    D^+(\Sigma)= \{ p \in M \, ; \, & \text{every inextendible past directed} \\ & \text{causal curve through $p$ intersects } \Sigma \}.
   \end{align*}
 \end{enumerate}
 Moreover, if $(M,C)$ is proper then $M$ is smoothly diffeomorphic to $\RR \times \Sigma$ (the projection to $\RR$ is a completely uniform temporal function), all Cauchy hypersurfaces are diffeomorphic to $\Sigma$ and the fibers of the smooth projection to $\Sigma$ are smooth timelike curves.
\end{thm}

 The proof of Theorem~\ref{thm:ghtime} for closed and proper cone structures is 
 based on the earlier works of Fathi and Siconolfi~\cite{FaSi}*{Theorem 1.3} for (i)$\Longleftrightarrow$(ii) and of Bernard and Suhr~\cite{BeSu:time}*{Theorem 3, Corollary 1.8} for (i)$\Longleftrightarrow$(iii) and the splitting. Note that in (iii) we use the recent notation of Bernard and Suhr \cite{BeSu:gh}*{Definition 1.2} rather than calling it $h$-steep as in \cite{Min:cone}*{Section 2.2}.

\medskip
 Being a closed cone structure does, of course, not imply that the causal relation is closed. But global hyperbolicity and properness suffices.
 
\begin{lem}[\cite{Min:cone}*{Lemma 2.5, Proposition 2.19}]\label{lem:compactclosed}
 Let $(M,C)$ be a proper cone structure. Then compactness of causal diamonds implies that the causal relation is closed or, equivalently, that the sets $J^\pm(p)$ are closed for all $p \in M$.
\end{lem}

 This result thus shows that global hyperbolicity is stronger than $K$-causality (because $K^+ = J^+$, which is antisymmetric if causal) in the framework of proper cone structures independent of the use of time functions. It is also clear that other notions on the causal ladder can be defined equally well for cone structures, however, none of those notions are relevant for us.

\subsection{The null distance for cone structures}\label{ssec:conenull}

 Equipped with a causal structure and time function we are finally in a position to introduce a cone version of the null distance of Sormani and Vega~\cite{SV}. The definitions of piecewise causal paths, null lengths, and null distance of Section \ref{ssec:nulld} carry over verbatim to closed cone structures. For the corresponding results to hold, however, the properness of cone structures is crucial.
  
\begin{lem}\label{lem:conn}
 Let $(M,C)$ be proper cone structure. Then there is a piecewise timelike (causal) curve between any two points of $M$.
\end{lem}

\begin{proof}
 The proof is analogous to that of \cite{SV}*{Lemma 3.5} and relies on the fact the the chronological future/past sets are open (and nonempty) by Proposition~\ref{prop:IJ}. For the sake of completeness we recall the full argument.
 
 For each $x \in M$ the sets $I^\pm(x)$ are open (and nonempty) by Proposition~\ref{prop:IJ} and hence $\{ I^-(x) \, ; \, x \in M \}$ is an open cover of $M$. Since $M$ is connected there is a continuous path $\alpha \colon [0,1] \to M$ between any to points $p,q \in M$. Because $\alpha([0,1])$ is compact it can be covered by finitely many sets $I^-(x_{2i})$, $i=0,1,\ldots,m$. Without loss generality we may assume that $p \in I^-(x_0)$, $q \in I^-(x_{2m})$ and $I^-(x_{2i}) \cap I^-(x_{2i+2}) \neq \emptyset$. Since $p \in I^-(x_0)$ there is a future directed timelike curve $\beta_0$ from $p$ to $x_0$. For all $0 \leq i<m$ fix a point $x_{2i+1} \in I^-(x_{2i}) \cap I^-(x_{2i+2}) \neq \emptyset$. Then there is a past directed timelike curve $\beta_{2i+1}$ from $x_{2i+1}$ to $x_{2i}$ and a future directed timelike curve $\beta_{2i+2}$ from $x_{2i+1}$ to $x_{2i+2}$. The concatenation $\beta = \beta_0 \cdots \beta_{2m}$ is a piecewise timelike curve from $p$ to $q$. 
\end{proof}

 Lemma~\ref{lem:conn} immediately implies that $\hat d_\tau$ is a pseudometric on $M$ but even in Minkowski space more is needed. The following result is the cone version of Theorem \ref{thm:defnulldist}.
 
\begin{thm}\label{thm:dttop}
 Let $(M,C)$ be a proper cone structure with a locally anti-Lipschitz time function $\tau \colon M \to \RR$. Then the null distance $\hat d_\tau$ is a metric on $M$ that induces the manifold topology.
\end{thm}

\begin{proof}
 By Lemma \ref{lem:conn} $\hat d_\tau \colon M \times M \to [0,\infty)$, and it is clear that $\hat d_\tau$ is symmetric and satisfies the triangle inequality. It thus remains to prove the distinguishing property, and that $\hat d_\tau$ induces the manifold topology.
 
 Let $x \in M$, $h$ a Riemannian metric, and $U$ the neighborhood of $x$ so that the anti-Lipschitz condition
 \[
  p \leq q \Longrightarrow d_h(p,q) \leq \tau(q)-\tau(p)
 \]
 holds for all $p,q \in U$. Suppose $y \neq x$. Let $V \subseteq U$ be a relatively compact open neighborhood of $x$ such that $y \not\in \overline{V}$. By Lemma~\ref{lem:conn} there is a piecewise causal path $\beta \colon [0,1] \to M$ from $x$ to $y$. Let $z = \beta(s_0) \in \partial V$ be the first point on $\beta$ that meets $\partial V$. Then
 \[
  \hat L_\tau(\beta) \geq \hat L_\tau(\beta|_{[0,s_0]}) \geq \hat d_\tau(p,z),
 \]
 and therefore $\hat d_\tau (x,y) = \inf \hat L_\tau(\beta) \geq d_\tau(x,\partial V) > 0$. Hence $\hat d_\tau$ is definite.
 
 It remains to be shown that $\hat d_\tau$ induces the manifold topology. As in \cite{SV}*{Proposition 3.14} one can show that the continuity of $\tau$ naturally implies the continuity of $\hat d_\tau$. That this proof carries over rests on the crucial fact that for proper cone structures $(\inte C)_p$ is nonempty for all $p \in M$ and that for every vector $v \in (\inte C)_p$ there is a timelike curve passing through $p$ with velocity $v$ (see the proof of Proposition \ref{prop:IJ}). Hence the topology induced by $\hat d_\tau$ is coarser than the manifold topology. As in \cite{SV}*{Proposition 3.15} it follows that it is also finer.
\end{proof}

 Theorem~\ref{thm:dttop} furthermore implies, as already established in \cite{AB}*{Theorem 1.1, Proposition 3.8} in the Lorentzian setting, that the length structure respects the manifold topology and that the null distance is an intrinsic metric.

\subsection{A completeness-compactness theorem for cone structures}
\label{ssec:coneHR}

 As in Theorem \ref{ABthm} it follows that the null distances of ``complete'' time functions are complete. Time functions that satisfy this condition exist for globally hyperbolic cone structures by Theorem \ref{thm:ghtime}. Because both results are valid for proper cone structures we obtain a natural extension to a refined Lorentzian completeness-compactness result of Garc\'ia-Heveling and the author~\cite{BGH2}*{Theorem 4.2}. Theorem \ref{conethm} is a simplified version of the following result.
 
\begin{thm}\label{thm:gghcomplete}
 Let $(M,C)$ be a proper cone structure.
 \begin{enumerate}
  \item Suppose $\tau$ is a time function such that the metric space $(M,\hat d_\tau)$ is complete. Then $\tau$ is a Cauchy time function and $(M,C)$ is globally hyperbolic.
  \item Suppose $(M,C)$ is globally hyperbolic. Then there exists a completely uniform (weak) temporal function $\tau \colon M\to\RR$, and for all such time function, the corresponding metric space $(M,\hat d_\tau)$ is complete. 
 \end{enumerate}
\end{thm}

\begin{proof}
 (i) It remains to be shown that $\tau$ is a Cauchy function, then global hyperbolicity follows from Theorem~\ref{thm:ghtime} (ii)$\Longrightarrow$(i). Suppose $\tau$ is not a Cauchy time function. Then there exists, without loss of generality, an inextendible future directed causal curve $\gamma \colon \RR \to M$ such that $A := \sup (\tau \circ\gamma) < \infty$. Consider the sequence of points $p_n = \gamma(n)$. For any $n,m \in \NN$
 \[
  \hat d_\tau(p_n,p_m) = |\tau(p_m) - \tau(p_n)|.
 \]
 Since the sequence $(\tau(p_n))_n$ is strictly increasing it converges to $A$ and therefore is a Cauchy sequence in $\RR$. In other words, for every $\varepsilon >0$ there is an $N \in \NN$ such that for all $n,m \geq N$ we have
 \[
  \hat d_\tau(p_n,p_m) = |\tau(p_m)-\tau(p_m)| < \varepsilon,
 \]
 thus $(p_n)_n$ itself is a Cauchy sequence in $(M,\hat d_\tau)$ and therefore most converge to a point $p \in M$ by completeness. Thus $\gamma$ is future extendible, a contradiction to the assumption. Hence $\tau$ must be a Cauchy time function.

 (ii) By Theorem~\ref{thm:ghtime} a completely uniform temporal function exists. As in \cite{BGH2}*{Theorem 4.2 (ii)} it follows from the cone version of Theorem \ref{ABthm} that $(M,\hat d_\tau)$ is complete.
\end{proof}

\section{Semi-Riemannian spacetimes}
\label{sec:sR}

 Thanks to Theorem~\ref{prop:equivLR} we have seen that Theorem~\ref{BGHthm} is a conformal generalization of the metric Hopf--Rinow Theorem~\ref{HRthm} (b)$\Longleftrightarrow$(c) to Lorentzian signature. A natural question that arises is whether one can iterate this procedure and also obtain a semi-Riemannian version of Theorem~\ref{BGHthm}. In this section we show that the answer is yes. To this end we introduce semi-Riemannian spacetimes and initiate their study in Section~\ref{ssec:sRspacetimes}. We provide several important examples and see that their existence is intimately tied to a deep and largely open problem in differential/algebraic topology. In Section \ref{ssec:sRcone} we answer the above question by showing that semi-Riemannian spacetimes admit continuous proper cone structures that are conformally invariant. Thus Theorem \ref{conethm} yields a true conformal and semi-Riemannian version of the metric Hopf--Rinow Theorem, namely Corollary \ref{semiRcor}. In Section \ref{ssec:sRGH} we further investigate if and how closely (stable) causality as well as global hyperbolicity are related for products (with regards to dimension and signature) and find, amongst others, that some but not all parts of Theorem~\ref{prop:equivLR} (i) can be recovered.

\subsection{Definition, existence, and examples}\label{ssec:sRspacetimes}

 Recall from the introduction that we consider the following class of semi-Riemannian manifolds $(M,g)$ with signature $(n-\nu,\nu)$, which interpolates between all manifolds equipped with a metric of positive definite signature $(n,0)$ and parallelizable manifolds equipped with a metric of negative definite signature $(0,n)$.

\begin{defn}\label{def:sRspacetime}
 Let $(M,g)$ be a semi-Riemannian manifold with constant index $0 \leq \nu \leq n = 
 \dim M$. We say that $M$ is \emph{time frame 
 orientable} if it admits $\nu$ continuous vector fields $X_i \in \vf(M)$ that in each $p \in M$ span a negative definite $\nu$-dimensional subspace of $T_pM$.
 
 If $(M,g)$ is time frame orientable and equipped with a fixed set $X = \{X_i \, ; \, i=1,\ldots,\nu\}$ of such vector fields, we say that it is \emph{time frame oriented} and call $(M,g,X)$ a \emph{semi-Riemannian spacetime} or, more specifically, a \emph{$(n-\nu,\nu)$-spacetime}.
\end{defn}

 In this section we discuss a vast list of examples and nonexamples of semi-Riemannian spacetimes and obtain some statements about the existence of such semi-Riemannian structures. We will see that there is a novel and interesting link of semi-Riemannian geometry to a challenging problem studied in differential and algebraic topology.
 
 \medskip
 By definition, any Riemannian manifold is time orientable and a $(n,0)$-spacetime with $X = \emptyset$. Every Lorentzian manifold that can be time oriented by a timelike vector field is a $(n-1,1)$-spacetime and a spacetime in the usual Lorentzian sense (cf.\ Section \ref{sec:Lbasics}). The use of the word "spacetime" in Definition \ref{def:sRspacetime} is fully justified in Section \ref{ssec:sRcone} where we see that the cone structure imposed by $X$ coincides with the usual one.
 It is well-known that not every Lorentzian manifold is time orientable, and that not every manifold even admits a Lorentzian metric. Using homotopy theory Steenrod first characterized the existence of a semi-Riemannian metric on a given (closed) manifold by the existence of a corresponding continuous tangent subbundle. We recall his argument in modern terminology. 
 
\begin{thm}[\cite{Ste}*{Theorem 40.11}]\label{Steenrodthm} 
  Let $M$ be a smooth manifold. The following are equivalent:
  \begin{enumerate}
   \item $M$ admits a semi-Riemannian metric of index $\nu$,
   \item $TM$ admits a subbundle of rank $\nu$.
  \end{enumerate}
\end{thm}

\begin{proof}
 Let $n = \dim M$ and $0 \leq \nu \leq n$.
 
 (i)$\Longrightarrow$(ii) Suppose $g$ is a semi-Riemannian metric of index $\nu$. This means that the tangent bundle $TM$ is associated to a principal $O(n - \nu, \nu)$ bundle by the standard representation of the group $O(n-\nu,\nu)$. The subgroup $O(n-\nu)\times O(\nu)$ is a maximal compact subgroup, and hence a deformation retract, of $O(n - \nu, \nu)$. So $TM$ can be associated to the product of the standard representation of $O(n-\nu)\times O(\nu)$. This implies that there are tangent subbundles $\xi$ and $\eta$  of ranks $\nu$ and $n - \nu$, respectively, satisfying
 $TM = \xi \oplus \eta$. 
 
 (ii)$\Longrightarrow$(i) Suppose $\xi$ is a rank-$\nu$ subbundle of $TM$. Let $h$ be any Riemannian metric $M$. Since $M$ is paracompact
 \[
  \xi^\perp := \bigsqcup_{p \in M} \{p\} \times \xi_p^\perp
 \]
 is a subbundle of rank $n-\nu$ of $TM$ via the natural projection $\xi^\perp \to M$ and such that 
 \[
  TM = \xi \oplus \xi^\perp. 
 \]
 A semi-Riemannian metric of index $\nu$ is then given by setting
 \[g = -h|_\xi \oplus h|_{\xi^\perp} = h - 2 h|_{\xi \times \xi}. \qedhere\]
\end{proof}
 
 In general, the frame bundle associated to the tangent subbundle (ii) does \emph{not} admit a global section, and thus the semi-Riemannian metric obtained in Theorem \ref{Steenrodthm} may \emph{not} be time frame orientable. Interestingly, if $\nu = 1$ the \emph{existence} of a Lorentzian metric implies the \emph{existence} of a time oriented Lorentzian metric (the converse being trivially true). After having discussed topological constraints characterizing such an existence we show in Example \ref{ex:s2} that the existence of a semi-Riemannian metric of index $\nu \geq 2$ on $M$ does \emph{not} imply the existence of a $(n-\nu,\nu)$-spacetime structure on $M$. We can, however, independently characterize the existence of a spacetime metric by the existence of a suitable global partial frame. Note that the condition of admitting such structures is invariant under $C^1$ diffeomorphisms.

\begin{thm}\label{thm:exst}
 Let $M$ be a smooth manifold of dimension $n$. The following are equivalent:
 \begin{enumerate}
  \item $M$ admits a $(n-\nu,\nu)$-spacetime metric,
  \item $M$ is \emph{parallelizable of degree $\nu$}, i.e., there exist $\nu$ everywhere linearly independent continuous vector fields on $M$ (called a \emph{tangent $\nu$-frame}).
 \end{enumerate}
\end{thm}

\begin{proof}
 (i)$\Longrightarrow$(ii) is trivial.
 
 (ii)$\Longrightarrow$(i)
 Given linearly independent vector fields $X_1,\ldots,X_\nu$, consider the rank-$\nu$ subbundle of $TM$ given by $\xi_p = \spann (X_1(p),\ldots,X_\nu(p))$. Follow the proof of Theorem \ref{Steenrodthm} to construct the metric $g$.
\end{proof}

\begin{cor}
 Every $n$-dimensional parallelizable manifold admits a $(n-\nu,\nu)$-spacetime structure for all $0 \leq \nu \leq n$. \hfill \qed
\end{cor}

\begin{ex}[Parallelizable manifolds]\label{ex:par}
 The following classes of manifolds are well-known to be parallelizable and thus can be turned into $(n-\nu,\nu)$-spacetimes for any $\nu$:
 \begin{enumerate}
  \item $\RR^n$, $\CC^n$, and all other finite-dimensional vector spaces over $\RR$ (see Section \ref{ssec:Rnu} for more properties),
  \item Lie groups (by choosing a basis at the identity and using group translations to move it around),
  \item closed orientable $3$-manifolds (by a result of Stiefel \cite{Sti:parallel}*{Satz 21}, see also the recent ``bare hands'' proof in \cite{BL:3mf}),
  \item $\mathbb{S}^1$, $\mathbb{S}^3$ (because $\mathbb{S}^3=\mathrm{SU}(2)$ is a Lie group), $\mathbb{S}^7$ (shown independently by Hirzebruch, Kervaire, and by Bott and Milnor in 1958), and no other spheres (see Adam's Theorem \ref{thm:spherevfs} below),
  \item closed $\pi$-manifolds (manifolds with trivial normal bundle when embedded in high dimensional Euclidean space, a concept due to Whitehead) of dimension $n$ are either parallelizable or have the same maximal $\nu$ as the corresponding $\mathbb{S}^n$ (this and related results are collected in \cite{Tho:vfs}*{Theorem 11}),
  \item products of parallelizable manifolds.
 \end{enumerate}
\end{ex}

 Parallelizable manifolds are useful for obtaining general semi-Riemannian spacetimes via the following (warped) product construction.

\begin{ex}[Semi-Riemannian warped products]\label{srprod}
 Let $M$ be an $n$-dimen\-sional manifold equipped with a $(n-\nu,\nu)$-spacetime structure. Suppose that $\Sigma$ is an $m$-dimensional manifold, equipped with an arbitrary Riemannian metric $\sigma$, and let $f \colon M \to (0,\infty)$ be a smooth function. Then by Theorem \ref{thm:exst} the warped product
 \[
  M \times_f \Sigma = (M \times \Sigma, g + f^2 \sigma)
 \]
 is a $(m+n-\nu,\nu)$-spacetime with the same time frame orienting vector fields (modulo pullback along $\pi_M$).
 Similarly, the warped product of $M$ with another $(m-\rho,\rho)$-spacetime $(N,k)$ is also a $(m+n-\rho-\nu,\rho+\nu)$-spacetime. Of course, this construction can be iterated. 
\end{ex}

 Let us return to the question of existence of a $(n-\nu,\nu)$-spacetime structure on a given manifold $M$. The corresponding topological problem, i.e., that of characterizing the existence of a tangent $\nu$-frame in condition (ii) of Theorem \ref{thm:exst}, is well-known and still largely open since almost 100 years. We shall provide some necessary conditions and discuss some progress that has been made since, and recall some of the important notions that have been developed in this context.
 
 \medskip
 The most basic necessary condition for the existence of a tangent $\nu$-frame is an immediate consequence of the Poincar\'e--Hopf Theorem \cite{Ho:vf} (an extension of the Hairy Ball Theorem of Poincar\'e \cite{Po:S2} for $\mathbb{S}^2$ from 1885 and of Brouwer \cite{Br: Sn} for $\mathbb{S}^{2n}$ from 1912) 
 since we certainly demand continuity and nowhere vanishing of at least one vector field. The full characterization in the $\nu=1$ case was shown in Markus \cite{Mar}*{Theorem 3}. 
 
\begin{thm}\label{thm:exSW}
 Let $M$ be a closed manifold that admits a $(n-\nu,\nu)$-spacetime structure for $\nu \geq 1$. Then $\chi(M) =0$. If $\nu = 1$ then this condition is also sufficient.
\end{thm}

 If $\nu=1$ one can furthermore show that the existence of Lorentzian metric on a given manifold implies the existence of a time oriented Lorentzian metric by moving to a double covering and applying Theorem \ref{thm:exSW} (see, for instance, \cite{ON}*{p.\ 149}). This is not the case for $\nu \geq 2$ as the following simple example demonstrates (unless, for instance, $M=\mathbb{S}^n$ and $2\nu \leq n$ \cite{Ste}*{Theorem 27.16}).

\begin{ex}[Manifold admitting a semi-Riemannian metric but no spacetime structure]\label{ex:s2}
 Consider $M = \mathbb{S}^2 \times \mathbb{S}^2$ and let $\sigma$ denote the standard sphere metric on $\mathbb{S}^2$. Then $g = - \sigma \oplus \sigma$ is clearly a semi-Riemannian metric of signature $(2,2)$. Since the Euler characteristic is $\chi(\mathbb{S}^2 \times \mathbb{S}^2) = \chi(\mathbb{S}^2)^2 = 4$ it follows from Theorem \ref{thm:exSW} that $M$ does not admit any semi-Riemannian spacetime structure of index $\nu=1,2,3,4$.
\end{ex}

 The refined notion of parallelizability of degree $\nu$ used in Theorem \ref{thm:exst} (ii) already appears in the seminal thesis of E.\ Stiefel \cite{Sti:parallel} 
 from 1935, supervised by H.\ Hopf, which he opens precisely by posing the following question in the closed case:
 \begin{quote}
 When does an $n$-dimensional manifold $M$ admit a tangent $\nu$-frame for $1 \leq \nu \leq n$?
 \end{quote}
 Neither Stiefel nor followers did so far succeed in fully characterizing the existence of a tangent $\nu$-frame for $\nu >1$. But already Stiefel obtained landmark results on the existence of tangent $\nu$-frames with certain types of singularities (points where the vector fields become linearly dependent or discontinuous) and introduced a sequence of obstruction classes in cohomology (independently and shortly afterwards Whitney \cite{Whi:top} studied the analogous classes for sphere bundles, the modern axiomatic definition of the Stiefel--Whitney classes for vector bundles is due to Hirzebruch). Although we are not interested in understanding the structure of singularities in this work, Stiefel's and subsequent results could provide a very fruitful starting point to pursue in the future. The reason for this is that we already know that understanding degeneracies in Lorentzian manifolds is related to finding an admissible notion for topology change of acausal slices, a topic that is highly relevant for quantum gravity and nonsmooth Lorentzian geometry (see, for instance, \cites{BDGS,BGH1,GH:top,Hor}).

 One of the most basic results for the existence of a tangent $\nu$-frame without singularities can in modern terminology and thanks to Theorem \ref{thm:exst} be formulated as follows.

\begin{thm}[Stiefel \cite{Sti:parallel}*{Satz A'$_m$}]\label{cor:wn}
 If a smooth manifold $M$ admits a $(n-\nu,\nu)$-spacetime structure then the top $\nu$ Stiefel--Whitney classes of the tangent bundle vanish, i.e., $w_{n-\nu + 1}(TM) = \ldots = w_{n}(TM) = 0$.
\end{thm}

 The necessary condition of Theorem~\ref{cor:wn} is very weak and far from being sufficient. For instance, any unit sphere $\mathbb{S}^n$ satisfies $w_1(T\mathbb{S}^n)=\ldots=w_n(T\mathbb{S}^n)=0$ because $w_0(T\mathbb{S}^n)=1$ and the total Stiefel--Whitney class is $w(T\mathbb{S}^n)=1$ \cite{Hus}*{Chapter 16, Proposition 4.4}, but we have already mentioned in Example \ref{ex:par} that only spheres of dimension $n=1,3,7$ are parallelizable.
 In 1962 Adams managed to settle the question on the maximal degree of parallelizability for spheres, based on earlier necessary criteria due to the Hurwitz--Radon--Eckman Theorem in linear algebra and James.
 
 \begin{thm}[Adams \cite{Ada:spherevfs} and James \cite{Jam:spherevfs}]\label{thm:spherevfs}
  The $n$-dimensional unit sphere $\mathbb{S}^n \subseteq \RR^{n+1}$ admits exactly $\nu(n) = 2^c + 8d -1$ linearly independent vector fields, where $c$ and $d$ are given implicitly by $n+1 = (2a +1)2^b$ and $b = c + 4d$ for $a,b,c,d \in \mathbb{Z}$, $0 \leq c  \leq 3$.
 \end{thm}
 
 Much less is known for general manifolds. Let us briefly recall the progress that has been made beyond Theorem \ref{cor:wn} since the 1930s. For certain classes of manifolds with particular dimensional restrictions $n$ (often $\operatorname{mod} 4$, spin, etc.) and for small $\nu$ knowledge about characteristic classes provides important necessary conditions. Thomas \cite{Tho:vfs} reviews many classical results of Atiyah, Bott, Frank, Hirzebruch, Hopf, Kervaire, Mayer, Milnor, Steenrod, Whitehead, Wu and many others from the 1950s and 1960s and discusses topological invariants that govern whether one can turn a tangent $\nu$-field with finite singularities into a regular one without singularities and collects restrictions and obstructions for the existence of a tangent $2$-field with singularities. In the 1970s Atiyah and Dupont \cites{AtDu,Du:K} from the 1970s, using the Atiyah--Singer Index Theorem (see also \cite{At}), obtain particularly strong results for the existence of tangent $2$- and $3$-fields for closed oriented manifolds allowing finite singularities. See also Crabb and Steer \cite{CrSt} using the same technique and Koschorke \cite{Ko:vflecture} treating these cases for nonorientable closed manifolds using a different approach. The more recent paper by B\"okstedt--Dupont--Svane \cite{BDS:cobord} gives a nice overview of known results and for $\nu=4,5,6,7$ (under additional assumptions on $M$ and $n$, such as spin) computes the index, whose vanishing characterizes the existence of a $\nu$-frame, and thereby shows that it is a global invariant.
 
\medskip
 Despite the huge amount of fascinating and deep mathematical work that has been carried out on the existence of tangent $\nu$-frames and singularities in the \emph{closed} manifold setting the equivalent problem for and also the topological structure of open (noncompact without boundary) manifolds is much less understood. For us the noncompact situation is significantly more important and in the hope to spark some interest in this problem we collect some known constructions and basic ideas that could be useful for analyzing the noncompact problem. 

\begin{rem}[Open manifolds]
 From a semi-Riemannian perspective the situation is better if compactness is dropped. On an open manifold there \emph{always} exists a nowhere vanishing continuous vector field and thus a Lorentzian spacetime structure (mentioned already in \cite{Mar}) because one can simply ``sweep out'' the zeroes of a vector field with isolated zeroes to infinity, for instance, by means of a suitable compact exhaustion and induction. It is not immediately obvious, though likely, that a similar procedure can be applied to obtain several linearly independent continuous vector fields. One may also be able to employ the homotopy equivalence property of deformation retractions: Smooth manifolds have the homotopy type of CW complexes. Whitehead showed that for open $n$-dimensional manifolds $M$ there is a subcomplex of dimension $k \leq n-1$ onto which $M$ deformation retracts. By Theorem \ref{cor:wn} this is a necessary condition for having $n-k$ linearly independent vector fields. 
 For more recent and sophisticated topological tools for analyzing noncompact manifolds see recent expository articles, for instance, \cite{Gui}.
 
 Masushita \cite{Ma:proc}*{p.\ 119} argues that Steenrod focused on the closed case because no conditions are needed for the existence of semi-Riemannian metrics of arbitrary index $\nu \geq 1$ in the open case. We could not find a proof for this claim in the literature, but even if we assume it it remains unclear whether the existence of a $(n-\nu,\nu)$-spacetime structure would also follow (in the closed case it does not because of Example \ref{ex:s2}). 
 
 If additional conditions for the existence of a tangent $\nu$-frame on an open manifold are needed, our setup could ultimately also shed some light on them. This is because the causal structure (which we explore in Sections \ref{ssec:sRcone} and \ref{ssec:sRGH}) and topology for $(n-\nu,\nu)$-spacetimes are beautifully intertwined and are able to capture global features of a manifold. We will, for instance, see in Section \ref{ssec:sRcone} that the results of Section \ref{sec:cone} can be applied to $(n-\nu,\nu)$-spacetimes. Corollary \ref{cor:closedcurves} then implies that there \emph{always} exist closed timelike (thus directly related to the chosen tangent $\nu$-frame!) curves on closed manifolds that admit a semi-Riemannian spacetime structure.
\end{rem}

\subsection{Semi-Riemannian spacetimes admit conformal proper cone structures}\label{ssec:sRcone}

 Given a spacetime structure on a manifold there are several ways to distinguish between future and past, and study causality. We pick a natural way to do so which for $\nu =1$ agrees with the usual Lorentzian one and show that it leads to proper cone structures.
 
\begin{defn}\label{def:future}
 Let $(M,g)$ be a $(n-\nu,\nu)$-spacetime with time frame orientation given by the vector fields $X_1,\ldots,X_\nu$. A vector $v \in T_pM \setminus \{0\}$ is said to be
 \begin{enumerate}
  \item \emph{future directed precausal} if 
 \begin{align}\label{ineqs:g}
   g_p(v,v) \leq 0, \quad \text{and} \quad g_p(v,X_i(p)) \leq 0, \,\, i=1,\ldots,\nu,
 \end{align}
  and \emph{future directed causal} if it is in the closed convex hull of the set of future directed precausal vectors,
  \item \emph{future directed pretimelike} if 
 \begin{align}
   g_p(v,v) < 0, \quad \text{and} \quad g_p(v,X_i(p)) < 0,  \,\, i=1,\ldots,\nu
 \end{align}
 and \emph{future directed timelike} if it is in the (open) convex hull of thet set of future directed pretimelike vectors.
 \end{enumerate}
 Similarly, we call a vector $v \in T_pM$ \emph{past directed} \emph{(pre)causal/(pre)timelike} if $-v$ is future directed (pre)causal/(pre)timelike.
\end{defn}


\begin{lem}\label{lem:onneg}
 Let $(M,g)$ be a $(n-\nu,\nu)$-spacetime with time frame orientation given by the vector fields $X_1,\ldots,X_\nu$.  If $v \in T_pM \setminus \{0\}$ is a future directed causal vector then there exists at least one $j \in \{1,\ldots,\nu\}$ such that
 \[
  g_p(v,X_j(p)) < 0.
 \]
\end{lem}

\begin{proof}
 Suppose there is a vector $v \in T_pM \setminus \{0\}$ that satisfies
 \begin{align}\label{null}
  g_p(v,X_i(p))=0, \qquad \text{for \emph{all} } i=1,\ldots,\nu.                                                                               
 \end{align}
 The subspace $\xi_p = \spann (X_1,\ldots,X_\nu)$ is a negative definite subspace of $(M,g)$ and $T_pM = \xi_p \oplus \xi_p^\perp$ with $g|_{\xi_p^\perp}$ being positive definite. Thus we can write $v$ uniquely as a sum of vectors $u \in \xi_p$ and $w \in \xi_p^\perp$. Assumption \eqref{null} implies that $g(u,u)=0$, i.e., $u=0$. If $v$ is precausal it satisfies $g(v,v)\leq0$ and hence
 \begin{align*}
  0 \geq g_p(v,v) = g_p(u,u) + g_p(w,w) = g_p(w,w) \geq 0,
 \end{align*}
 i.e., also $w=0$, a contradiction to the assumption that $v\neq 0$. If $v$ is causal it is a convex combination $v = \sum_{l=1}^k \alpha_l v_l$ of precausal vectors $v_l$ with the same property, hence all $w_l = 0$, which implies $w=0$ and again contradicts $v \neq 0$. Thus there must be at least one index $j\in\{1,\ldots,\nu\}$ for which $g_p(v,X_i(p))<0$.
\end{proof}

\begin{rem}[Alternative definition]\label{def:altfuture}
 If $(M,g)$ is a semi-Riemannian spacetime then 
 for each $p \in M$ we can write
 \[
  T_pM = \xi_p \oplus \xi_p^\perp,
 \]
 where $\xi_p = \spann (X_1(p),\ldots,X_\nu(p))$. Thus we can write every tangent vector $v \in T_pM \setminus \{0\}$ as $v = u + w$ with unique $u \in \xi_p$ and $w \in \xi_p^\perp$. Suppose $u = \sum_{i=1}^\nu u^i X_i(p)$ is the unique representation in terms of the basis of $\xi_p$. Then we could also define $v$ to be \emph{future directed} if $u^i \geq 0$ for all $i=1,\ldots,\nu$. If the vector fields $X_i$ are orthogonal this definition agrees with Definition~\ref{def:future}. By the Gram--Schmidt orthogonalization procedure we can indeed always assume orthonormality of the vector fields $X_1,\ldots,X_\nu$ but we chose not to require this condition in our setup.
\end{rem}

\begin{rem}[Lorentzian vs.\ semi-Riemannian]\label{rem:unrelated}
 Unlike in the Lorentzian case, not every vector $v \neq 0$ satisfying $g_p(v,v) \leq 0$ is \emph{either} future \emph{or} past directed if $\nu > 1$. In fact, the set of ``undirected (pre)causal'' vectors is actually bigger than the directed ones. Note also that vectors that are causal but not precausal satisfy $g_p(v,v)>0$.
 Another way of seeing this is that in the Lorentzian setting a different choice of time orientation defining vector field at most exchanges the future with the past, while in the semi-Riemannian setting it generally leads to an entirely different causal structure and a particular choice of tangent $\nu$-frame is more intimately linked to the topology of $M$. See also Example \ref{ex:ghngh}.
\end{rem}

 Nonetheless, in analogy to the Lorentzian setting, we show that the definition of future/past directed causal vectors defines a proper cone structure on $M$. The regularity of the metric tensor and the vector fields controls the regularity of the cone structure and ensures that it is topologically well behaved. We restate Theorem \ref{thm:sRcones} of the introduction in a more refined way as follows.

\begin{thm}\label{thm:propercone}
 Let $(M,g)$ be a $(n-\nu,\nu)$-spacetime with $0<\nu\leq n$. The map
 \begin{align*}
  p \mapsto C_p := \{ v \in T_p M \setminus \{ 0\} \, ; \, v ~\text{is future directed causal} \}
 \end{align*}
 defines a continuous proper cone structure on $M$ in the sense of Definition~\ref{def:constructure} with
 \[ (\inte C)_p = \inte C_p = \{ v \in T_p M \setminus \{ 0\} \, ; \, v ~\text{is future directed timelike} \}.\]
\end{thm}

 Before we prove Theorem \ref{thm:propercone} a few remarks are in order. Let us first recall the definition of continuity of set-valued maps that is used in the statement (see, for instance, Aubin and Cellina~\cite{AuCe:setvalued}*{Chapter 1}).

\begin{defn}\label{def:setcont}
 Let $F \colon X \to Y$ be a set-valued map. Then we say that
 \begin{enumerate}
  \item $F$ is \emph{upper semicontinuous} at $x_0 \in X$ if for any neighborhood $V$ of $F(x_0)$, there exists a neighborhood $U$ of $x_0$ such that $F(U) \subseteq V$.
  \item $F$ is \emph{lower semicontinuous} if for any net of elements $x_\mu$ converging to $x_0$ and for any $y_0 \in F(x_0)$, there exists a sequence of elements $y_\mu \in F(x_\mu)$ that converges to $y_0$.
 \end{enumerate}
 The map $F$ is \emph{continuous} if it is both upper and lower semicontinuous.
\end{defn}

 For cone structures $p \mapsto C_p$ on a manifold $M$ these properties can be checked on coordinate patches (see \cite{Min:cone}*{p.\ 12}): At every $p \in M$ we have a local coordinate system $\{x^\alpha\}$ on a neighborhood $U$ of $p$. Via the local trivialization of the tangent bundle $TU$ we obtain a splitting $U \times \RR^n$ and can compare sets over different tangents spaces this way. The notion of upper and lower semicontinuity are then applied to the set-valued map
 \[ F(p) = [ C_p \cup \{0\}] \cap \mathbb{B}^n, \]
 where $\mathbb{B}^n$ denotes the closed unit ball of $\RR^n$ (instead one could also work with the sphere subbundle and compare the set values $C_p \cap \mathbb{S}^n$ using the Hausdorff distance on $\mathbb{S}^n$; one can also say when $F$ is locally Lipschitz then).

 \medskip
 Furthermore, recall that properness of a cone structure requires a nonempty interior bundle. We thus construct a vector in the interior of the cone explicitly by making use of the negative definite version of the following result.

\begin{lem}\label{lem:intC}
 Let $V$ be a finite-dimensional inner product space, and let $(b_1,\ldots,b_m)$ be a basis for $V$. Then there exists a vector $v \in V \setminus \{ 0 \}$ such that
 \begin{align}\label{eq:vbi}
  \langle v, b_i \rangle > 0, \qquad i =1,\ldots,m.
 \end{align}
\end{lem}

\begin{proof}
 We proceed by induction in $m$. If $m=1$, then we can simply pick $v = b_1$.
 
 Suppose we have shown the result already for the inner product spaces of dimension $m-1$ and $\tilde v$ is the vector satisfying \eqref{eq:vbi} for $\widetilde V =\spann (b_1,\ldots,b_{m-1})$. By the Gram--Schmidt orthogonalization procedure there exists an orthogonal basis $(e_1,\ldots,e_m)$ of $V$ such that
 \[
  \spann(b_1,\ldots,b_k) = \spann(e_1,\ldots,e_k), \qquad k = 1,\ldots,m.
 \]
 In particular, this holds for $k=m-1$, and therefore \[ \langle b_m, e_m \rangle \neq 0, \] since $b_m$ is a basis vector and therefore not in $\widetilde V = \{ e_m \}^\perp = \spann(b_1,\ldots,b_{m-1})$. Choose $a \in \RR$ such that
 \[
  \langle b_m, e_m \rangle a > - \langle \tilde v, b_m \rangle
 \]
 and define
 \[
  v = \tilde v + a e_m.
 \]
 Clearly $v \neq 0$. We show that $v$ satisfies \eqref{eq:vbi}. For any basis vector $b_k$ with $k = 1,\ldots,m-1$ we have $b_k \in \widetilde V  = \{ e_m \}^\perp$, and hence \eqref{eq:vbi} follows directly from the induction hypothesis since
 \[
  \langle v, b_k \rangle = \langle \tilde v, b_k \rangle > 0. 
 \]
 For $b_m$ we have by choice of $a$ also
 \begin{align*}
  \langle v, b_m \rangle &= \langle \tilde v, b_m \rangle + a \langle b_m,e_m \rangle > 0,
 \end{align*}
 which completes the inductive step from $m-1$ to $m$.
\end{proof}

\begin{proof}[Proof of Theorem~\ref{thm:propercone}]
 Suppose the time frame orientation of $(M,g)$ is given by the continuous vector fields $X_1,\ldots,X_\nu \in \vf(M)$.
 
 We first investigate the properties of the cones pointwise. Let $p \in M$. By Definition~\ref{def:future}, $C_p := \overline{\operatorname{co}}~ \mathcal{C}_p$ with $\overline{\operatorname{co}}$ the closed convex hull operator and
 \[
  \mathcal{C}_p = G_p \cap \left( \bigcap_{i=1}^\nu C^i_p \right),
 \]
 for the sets
  \[ G_p := \{ v \in T_pM \setminus \{0\} \, ; \, g_p(v,v) \leq 0 \}\]
 and
 \[
   C^i_p := (X^\flat_i)^{-1}_p(-\infty,0]= \{ v \in T_pM 
   \, ; \, g_p(v,X_i(p)) \leq 0 \}.
 \]
 Bilinearity of $g$ and Lemma~\ref{lem:onneg} imply that for each $p \in M$ the set $\mathcal{C}_p$ is a sharp cone (it is convex iff $\nu =1$ or $\nu=n$). Furthermore, $\mathcal{C}_p$ is closed in $T_pM \setminus \{0\}$ as finite intersection of closed sets in the subspace topology of $T_p M \setminus \{0\}$. Each cone $\mathcal{C}_p$ has nonempty interior because for each $p \in M$, by Lemma~\ref{lem:intC} (in the negative definite case), there exists $v \in \spann(X_1(p),\ldots,X_\nu(p)) \setminus \{0\} \subseteq T_pM$ such that
 \[
  g_p(v,X_i(p)) < 0, \qquad i = 1,\ldots,\nu.
 \]
 Since $g_p(v,v)<0$ we also know that $v$ is timelike (and not just causal) and therefore
 \[
  v \in \inte G_p \cap \left( \bigcap_{i=1}^\nu \inte C_p^i \right) = \inte \mathcal{C}_p \neq \emptyset.
 \]
 
 Thus each $C_p = \overline{\operatorname{co}}~ \mathcal{C}_p$ is a proper (convex) cone and it remains to be shown that the cone structure $C$ is closed and proper not only pointwise but as subbundle of $TM \setminus (TM)_0$. Our approach makes use of the multivalued map $p \mapsto \mathcal{C}_p$ for which we establish continuity similar to Minguzzi \cite{Min:cone}*{Proposition 2.4} in the Lorentzian case and then transfer this property to $p \mapsto C_p$: We first prove closedness of $\mathcal{C}$ which by \cite{AuCe:setvalued}*{Theorem 1.1.1} immediately also implies upper semicontinuity. In a second step, we show lower semicontinuity of $\mathcal{C}$ directly. Continuity extends to $C$ via conic combinations and we conclude properness using \cite{Min:cone}*{Proposition 2.5}.

 Suppose $p \in M$ and $U$ is a chart neighborhood of $p$. The cone bundle $\mathcal{C}$ is characterized by nonpositivity of the function
 \[
  f(q,w) := \max \{ g_{\alpha\beta}(q)w^\alpha w^\beta, g_{\alpha\beta}(q)X_1^\alpha(q)w^\beta, \ldots, g_{\alpha\beta}(q)X_\nu^\alpha(q)w^\beta \},
 \]
 which is continuous on the associated local trivialization $U \times \RR^n$ of the tangent bundle. The interior is characterized by negativity of $f$.
 
 Note that
 \[
  (\mathcal{C} \cup (TM)_0) \cap TU = \{(q,w) \, ; \, q \in U, w \in T_qM, ~\text{and}~ f(q,w) \leq 0 \}
 \]
 is closed in the topology of $TU$ by continuity of $f$. This yields that $\mathcal{C} \cup (TM)_0$ is closed in $TM$, which implies upper semicontinuity of $p \mapsto \mathcal{C}_p$ (let $F$ in \cite{AuCe:setvalued}*{Theorem 1.1.1} be the distribution of compact unit coordinate spheres $\mathbb{S}^n$).

 For lower semicontinuity, since $\overline{\inte \mathcal{C}_p} = \mathcal{C}_p \cup \{0\}$, it is sufficient to consider the set-valued map $p \mapsto \inte \mathcal{C}_p$. Suppose $(p,v) \in \inte \mathcal{C}_p $. Then there is an $\varepsilon < 0$ such that $f(p,v) < \varepsilon < 0$. If now $p_n \to p$, then there is an $N \in \NN$ such that for all $n \geq N$ we have $f(p_n,v) < 0$ as well. Hence $(p_n,v) \in \inte \mathcal{C}_{p_n}$. For the first $N-1$ elements $p_n$ of the sequence one can also find vectors $v_n \in \inte \mathcal{C}_{p_n}$ by Lemma~\ref{lem:intC}. Furthermore, $v_n$ converges to $v$ since they agree for all $n \geq N$, which ensures the lower semicontinuity of $p \mapsto \mathcal{C}_p$. 
 
 Having shown lower and upper semicontinuity for the multivalued map $p \mapsto \mathcal{C}_p$ we conclude that it is continuous. In particular, the map $p \mapsto \mathcal{C}_p \cap \mathbb{S}^n$ into the compact generating sets of the cones $C_p$ is continuous with respect to the Hausdorff distance. Since the conic hull operation itself is continuous for sharp cones \cite{Orl:cone}*{Proposition 2} (this result also holds for compact generating sets), the map $p \mapsto C_p$ itself is continuous. By \cite{Min:cone}*{Proposition 2.5} continuity implies that $C$ is a proper cone structure. By \cite{Min:cone}*{Proposition 2.6} continuity and properness furthermore imply $(\inte C)_p = \inte C_p$ for all $p \in M$.
\end{proof}

 Thanks to Theorem \ref{thm:propercone} we can apply the tools of the theory of cone structures from Section \ref{sec:cone} to semi-Riemannian spacetimes using $p\mapsto C_p$. For instance, the causal and chronological relation behave as desired and we can use the most important tools of Lorentzian causality theory in the semi-Riemannian context\footnote{With Lipschitz regularity the (for $1 < \nu < n$) nonconvex precausal cone map $p\mapsto \mathcal{C}_p$ could also be used to obtain solutions to the corresponding differential inclusion but limit curves, which we implicitly use via \cite{Min:cone}, are generically only in the relaxed convexified differential inclusion \cite{AuCe:setvalued}*{Sections 2.3 and  2.4}.}. The steps on the causal ladder and, in particular, global hyperbolicity can be defined directly via the theory of cone structures and Theorems \ref{thm:ghtime} and \ref{conethm} apply. We obtain Corollary \ref{semiRcor} directly from Theorems \ref{thm:gghcomplete} and \ref{thm:exst}.
 
 That we actually gain information in our semi-Riemannian setup in contrast to general cone structures is evident from its conformally invariant setup.
 
 \begin{thm}\label{thm:confnulld}
  Let $(M,g)$ be a $(n-\nu,\nu)$-spacetime, $0 <\nu \leq n$, with future causal cone defined pointwise as in Definition \ref{def:future}. The notions of causal relation $J^+$ and the chronological relation $I^+$ (and the steps on the causal ladder) induced by the corresponding proper cone structure are conformally invariant, and so are time functions $\tau$ and the null distance $d_\tau$.
 \end{thm}
   
 \begin{rem}[$\nu = 0$]
  Although one can successful interpret Riemannian manifolds as $(n,0)$-spacetimes it is---as expected---not very insightful to study their causal structure. Since the cone $C = \emptyset$ is degenerate and thus all $J^+(p)=\{p\}$ it would mean that all Riemannian manifolds are globally hyperbolic and every continuous function is a time function (if one would generalize this definition at all, the theory of cone structures often disregards degenerate cones).
 \end{rem}
 
 Finally, note that the stability of global hyperbolicity (as shown for continuous proper cone structures by Fathi and Siconolfi \cite{FaSi}*{Theorem 1.2} and later extended by others, see also \cite{Min:cone}*{Theorem 2.39}) can in the semi-Riemannian context be interpreted not only in terms of a nearby metric tensor $\tilde g$ but also in terms of a nearby time frame orientation. We show a result via the second approach.

\begin{lem}
 Let $(M,g)$ be a globally hyperbolic $(n-\nu,\nu)$-spacetime, $0 < \nu \leq n$, with time frame orientation defining vector fields $X_1,\ldots,X_\nu$ and corresponding proper cone structure $C$. Then there exist time frame orientation defining vector fields $\widetilde X_1,\ldots,\widetilde X_\nu$ on $(M,g)$ such that the corresponding cone structure $\widetilde C$ satisfies $C \preccurlyeq \inte \widetilde C \cup (\partial \widetilde C \setminus \bigcup_{i=1}^\nu (\widetilde X_i^\flat)^{-1} (0))$ and $\widetilde C$ is also globally hyperbolic.
\end{lem}

\begin{proof}
 Let $X$ be the $g$-unital version of the vector field ${X_1 + \ldots + X_\nu}$. Then for any $\varepsilon \in (0,1)$ the vector fields
  \begin{align}\label{tXi}
   \widetilde{X_i} &= X_i + \varepsilon X, \quad i = 1,\ldots, \nu, 
  \end{align}
  define a time frame orientation such that for each future directed causal vector $v \in C_p$ in $M$ we have, by Lemma~\ref{lem:onneg}, that
  \begin{align*}
    g_p(v,\widetilde X_i(p)) < 0, \qquad i=1,\ldots,\nu.
  \end{align*}
 In particular, $v \in \inte \widetilde{C}_p \cup (\partial \widetilde{C}_p \setminus \bigcup_{i=1}^\nu (\widetilde X_i^\flat)_p^{-1} (0))$, where $\widetilde C$ denotes the cone structure corresponding to the time frame orientation $\widetilde X_1,\ldots,\widetilde X_\nu$.
 
 It remains to be shown that $\widetilde C$ can be adapted to be globally hyperbolic. Since $(M,C)$ is globally hyperbolic and globally hyperbolic cone structures are stable there exists a globally hyperbolic locally Lipschitz proper cone structure $C'$ such that $C \prec C'$ (see, for instance, \cite{Min:cone}*{Theorem 2.39}). Since $C'$, $C$ and $\widetilde C$ are all continuous cone structures we can continuously choose an $\varepsilon (p) >0$ on $M$ defining $\widetilde C$ as above by \eqref{tXi} such that $C \preccurlyeq \widetilde C \prec C'$. Since $C'$ is globally hyperbolic it admits a Cauchy time function $\tau$. Clearly, $\tau$ is also a Cauchy time function for $\widetilde C$ and hence it is globally hyperbolic by Theorem~\ref{thm:ghtime}.
\end{proof}

\subsection{Semi-Riemannian vector spaces}\label{ssec:Rnu}
 In this section we are concerned with the flat spacetime $\RR^{n-\nu,\nu}$. We first fix the canonical $(n-\nu,\nu)$-spacetime structure on $\RR^n$ (and $\CC^n$) by using the standard orthonormal global coordinate frame $(E_i)_i$.

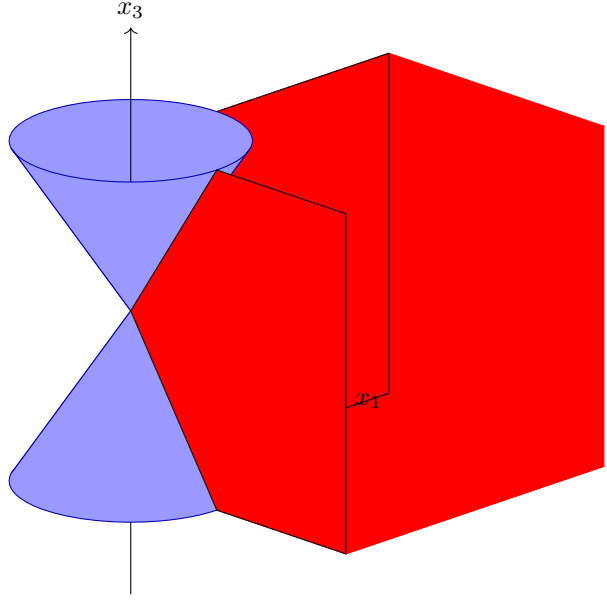
\begin{figure}[h]
 \begin{subfigure}[b]{0.49\textwidth}
   \centering
   \resizebox{\linewidth}{!}{
 \begin{tikzpicture}[tdplot_main_coords,scale=0.8]
  \coordinate (O) at (0,0,0);
  \draw[->] (0,0,0) -- (0,6,0) node[above right] {$x_2$};
  \fill[canvas is yz plane at x=5,fill=red,draw=red,fill opacity=0.4] (0,-3) -- (6,-3) -- (6,3) -- (0,3) -- (0,0) -- cycle; 
  \fill[canvas is xy plane at z=3,fill=red,draw=red,fill opacity=0.4] (0,0) rectangle (5,6);
  \fill[canvas is xy plane at z=-3,fill=red,draw=red,fill opacity=0.1] (0,0) rectangle (5,6);
  \draw[canvas is yz plane at x=0,fill=red, fill opacity=0.6] (2,-3) -- (6,-3) -- (6,3) -- (2,3) -- (0,0) -- cycle;
  \draw (0,0,-5) -- (O);
  \conefront[surface]{-3}{2}{-10}
  \coneback[surface]{3}{2}{10}
  \draw[->] (O) -- (0,0,5) node[above] {$x_3$};
  \conefront[surface]{3}{2}{10}
      \draw[->] (0,0,0) -- (5,0,0) node[below right] {$x_1$};
  \draw[canvas is xz plane at y=0, fill=red,opacity=0.6] (2,-3) -- (5,-3) -- (5,3) -- (2,3) -- (0,0) -- cycle;
\end{tikzpicture}}
\caption{}
\end{subfigure}
\begin{subfigure}[b]{0.41\textwidth}
   \centering
   \resizebox{\linewidth}{!}{
   \begin{tikzpicture}[tdplot_main_coords,scale=0.8]
  \coordinate (O) at (0,0,0);
  \draw[->] (0,0,0) -- (0,6,0) node[above right] {$x_2$};
  \fill[canvas is yz plane at x=5,fill=green,draw=green,fill opacity=0.4] (0,-3) -- (6,-3) -- (6,3) -- (0,3) -- (0,0) -- cycle; 
  \fill[canvas is xy plane at z=3,fill=green,draw=green,fill opacity=0.4] (2,0) -- (5,0) -- (5,6) -- (0,6) -- (0,2) -- cycle;
  \fill[canvas is xy plane at z=-3,fill=green,draw=green,fill opacity=0.4] (2,0) -- (5,0) -- (5,6) -- (0,6) -- (0,2) -- cycle;
  \draw[canvas is yz plane at x=0,fill=green, fill opacity=0.6] (2,-3) -- (6,-3) -- (6,3) -- (2,3) -- (0,0) -- cycle;
  \filldraw[draw=black, fill=green, fill opacity=0.4] (0,0,0) -- (2,0,3) -- (0,2,3) -- cycle;
  \filldraw[draw=black, fill=green, fill opacity=0.4] (0,0,0) -- (2,0,-3) -- (0,2,-3) -- cycle;
   \draw (0,0,-5) -- (O);
  \draw[->] (O) -- (0,0,5) node[above] {$x_3$};
      \draw[->] (0,0,0) -- (5,0,0) node[below right] {$x_1$};
  \draw[canvas is xz plane at y=0, fill=green,opacity=0.6] (2,-3) -- (5,-3) -- (5,3) -- (2,3) -- (0,0) -- cycle;
\end{tikzpicture}
   }
   \caption{}
\end{subfigure}
\caption{\small (A) In $\RR^{1,2}$ the cones $\mathcal{C}_p \subseteq T_p \RR^{3} \cong \RR^{3}$ (in red) are delimited by the intersection of the half spaces $(E_i^\flat)^{-1}(-\infty,0] = \{ v \in \RR^3 ; \, \langle v, E_i \rangle \leq 0 \}$ for $i=1,2$ with the "null cone" $\partial G = \{v \in \RR^3 \setminus \{0\} ; \, \langle v,v \rangle = 0\}$ (in blue). (B) The cones $C_p$ (in green) are the closed convex hulls of the cones $\mathcal{C}_p$. See also the proof of Theorem~\ref{thm:propercone} for the notation.}
\label{fig:cone}
\end{figure}

\begin{ex}[$\RR^{n-\nu,\nu}$]\label{rmn:spacetime}
 We equip the manifold $\RR^n$ with the standard scalar product of index $\nu$: With respect to the standard global coordinate frame $(E_1,\ldots,E_n)$, and $v = v^i E_i$, $w = w^i E_i$, we have a semi-Riemannian metric given by
 \[
  \langle v,w \rangle = - \sum_{i=1}^\nu v^i w^i + \sum_{j=\nu+1}^{n} v^j w^j. 
 \]
 The vector fields $E_1,\ldots, E_\nu$ are negative unital and orthonormal, and thus define a time frame orientation on $\RR^{n-\nu,\nu}$. Note, however, that for $\nu \geq 2$ these vector fields are \emph{not} future directed timelike (only causal). A cone $C_p$ is depicted in Figure~\ref{fig:cone}.
\end{ex}

\begin{ex}[$\CC^n$ as the real $(n,n)$-spacetime $\RR^{n,n}$]
 The complex vector space $\CC^n \cong \RR^n \oplus i \RR^n$ with scalar product
 \[
  \langle w,z \rangle = - \sum_{i=1}^n \Im (w^i) \Im (z^i) + \sum_{i=1}^n \Re (w^i) \Re (z^i)
 \]
 is isomorphic to the $(n,n)$-spacetime $\RR^{n,n}$ and the spacetime structure respects the complex nature of $\CC^n$.
\end{ex}
 
 We discuss and derive several properties of $\RR^{n-\nu,\nu}$ by hand. All cases $0 < \nu \leq n$ are covered. If $\nu = n$ the empty sum convention $\sum_{j=n+1}^n = 0$ is used.
 
 \subsubsection*{$I^+$ and $J^+$ relations}
 Two points $p,q \in \RR^{n-\nu,\nu}$ are future/past directed timelike/cau\-sally/light\-like related if and only if the straight line connecting them has the same causal character, i.e., if the corresponding vector $q-p \in \RR^{n-\nu,\nu} \cong T_0 \RR^{n-\nu,\nu}$ is future/past directed timelike/causal/lightlike. The convex hull operation can be carried out explicitly to characterize the causal relation in this setting (convex combinations can be split up with the Minkowski inequality, which singles out the $1$-norm on $\spann(E_1,\ldots,E_\nu)$ as necessary and sufficient upper bound for convexity). For $p = \sum_{i=1}^n p^i E_i$ and $q = \sum_{i=1}^n q^i E_i$ 
 \begin{align*}
  p \leq q \Longleftrightarrow\, & q^i \geq p^i \text{ for all }i=1,\ldots,\nu, \text{ and } \\
  &\sum_{i=1}^\nu q^i-p^i \geq \sqrt{\sum_{j=\nu+1}^{n} (q^j - p^j)^2}.
 \end{align*}
 If $p\neq q$ then automatically at least one $q^i > p^i$ (see also Lemma~\ref{lem:onneg}). Analogous characterizations hold for timelike and lightlike with all $>$ and at least one $=$ compared to the above, respectively (recall that lightlike means $q-p \in \partial C$).

\subsubsection*{Canonical time function}

 The function
 \[
  T(p) := \sum_{i=1}^\nu p^i
 \]
 is a smooth time function because by the above 
 \begin{align*}
  p < q \Longrightarrow \, & q^i \geq p^i \text{ for all } i=1,\ldots,\nu, \text{ and }\\
  & q^k > p^k \text{ for at least one } k=1,\ldots, \nu, \\
  \Longrightarrow \, & T(q) > T(p).
 \end{align*}
 It is easy to see that $T$ is, in fact, completely uniform Cauchy temporal: Since $dT = \sum_{i=1}^\nu  E^\flat_i$, for any $v \in T_p \RR^{n-\nu,\nu}$ future directed causal $v^i \geq 0$, $i=1,\ldots,\nu$, and
 \begin{align*}
  dT(v) = \sum_{i=1}^\nu v^i \geq \sqrt{\sum_{i=\nu+1}^n (v^i)^2},  
 \end{align*}
 which implies that it is bounded below by the Euclidean norm of $v$ since
 \begin{align*}
 2 dT(v) \geq \sum_{i=1}^\nu v^i + \sqrt{\sum_{i=\nu+1}^n (v^i)^2}\geq \sqrt{\sum_{i=1}^n (v^i)^2} = \| v\|_2.
 \end{align*}

 If $\nu=1$ then $T$ is the canonical time function used for Minkowski space.

\subsubsection*{Global hyperbolicity}

 Since $T$ is a Cauchy time function it follows from Theorem \ref{thm:ghtime} (ii)$\Longrightarrow$(i) that $\RR^{n-\nu,\nu}$ is globally hyperbolic. In connection with Theorem~\ref{prop:equivLR} that we aim to (partly) generalize in Section \ref{ssec:sRGH} it is insightful to prove global hyperbolicity also with ``bare hands'' using the Heine--Borel property of Euclidean space $\RR^n$.
 
 For any $p \in \RR^{n-\nu,\nu}$
  \[
  J^+(p) = \left(\bigcap_{j=1}^{\nu} \{ x^j \geq p^j \}\right) \cap \left\{ \sum_{i=1}^
\nu x^i-p^i \geq \sqrt{\sum_{i=\nu+1}^{n} (x^i-p^i)^2} \right\},
 \]
 and thus the causal future and past sets $J^\pm(p)$ are closed as an intersection of closed sets.  Hence for any $p$ and $q$ the causal diamond $J^+(p) \cap J^-(q)$ is closed. In what follows we show that $J^+(p) \cap J^-(q)$ is a bounded subset of the Euclidean space $\RR^{n}$. Since
 \[
  J^+(p) \cap J^-(q) \subseteq \bigcap_{j=1}^{\nu} \{ q^j \geq x^j \geq p^j \}
 \]
 the set is bounded in the first $\nu$ coordinate directions by future and past directedness. If $\nu=n$ we have shown boundedness of $J^+(p) \cap J^-(q)$. If $\nu < n$ boundedness in the remaining $n-\nu$ directions follows from causality. For $x \in  J^+(p) \cap J^-(q)$ and $i=1,\ldots,\nu$
 \[
  |q^i - x^i| + |x^i - p^i| = q^i - p^i \leq \max_{1\leq j\leq \nu} q^j - p^j =: \frac{R}{\nu}
 \]
  and thus
 \[
  \sqrt{\sum_{i=\nu+1}^{n} (x^i-q^i)^2} + \sqrt{\sum_{i=\nu+1}^{n} (x^i-p^i)^2} \leq 2R,
 \]
 meaning that $x$, and thus $J^+(p) \cap J^-(q)$, is also bounded in the remaining $n-\nu$ coordinate directions. 
 
 Finally, we apply the Heine--Borel property of the Euclidean space and conclude that $J^+(p) \cap J^-(q)$ is compact.
 
\subsubsection*{Null distance}

 By Theorem \ref{thm:confnulld} the null distance $\hat d_T$ is a conformally invariant metric that induces the topology of $\RR^{n-\nu,\nu}$. Sormani and Vega \cite{SV}*{Proposition 3.4} have shown that Minkowski space can easily be equipped with a non-locally anti-Lipschitz time function so that the corresponding null distance is only a pseudometric. We show that $\tau = T^{2k+1}$, $k \in \NN$, are analogous examples for $\RR^{n-\nu,\nu}$ as long as $\nu < n$. Let $p = (0,\ldots,0)$ and $q = (0,\ldots,0,1)$. By the above $p \nleq q$. The paths $\gamma^\pm \colon [0,1]\to M$,
 \[
  \gamma^\pm(t) \to (\pm t,0,\ldots,t)
 \]
 are causal, with $\gamma^+$ being future directed and $\gamma^-$ being past directed. By rescaling and concatenating such future and past directed causal paths we obtain piecewise causal paths $\beta_j \colon [0,1] \to [0,\frac{1}{2j}] \times \{ 0 \} \times \ldots \times [0,1]$ from $p$ to $q$ with $2j$ pieces of $\tau$-height $(\frac{1}{2j})^{2k+1}$, i.e.,
 \[
  \hat L_\tau (\beta_j) = \frac{2j}{(2j)^{2k+1}} = \frac{1}{4^k j^{2k}}.
 \]
 Thus
 \[
  \hat d_\tau (p,q) \leq \lim_{j\to\infty} \hat L_\tau(\beta_j) = 0,
 \]
 which means that $\hat d_\tau$ is not distinguishing with respect to any $\tau = T^{2k+1}$.
 
 This example furthermore shows that the case $\nu= n$ is special, which foreshadows some results of the following section.

 \medskip
 After having discussed $\RR^{n-\nu,\nu}$ with the canonical time frame orientation induced by $E_1,\ldots,E_\nu$ we conclude with a simple modification thereof, demonstrating that the choice (and orientation) of vector fields is crucial for the global causal properties of a semi-Riemannian spacetime.
 
\begin{ex}[Relevance of the time frame orientation]\label{ex:ghngh}
 On $\RR^{n-\nu,\nu}$ the chosen time frame orientation does not affect the causal properties of the spacetime. It does, however, if we mildly modify $\RR^{0,2}$ to be $M := \RR^{0,2} \setminus \{ (x_1,x_2) \, ; \, x_1=x_2>0 \}$ with the standard negative definite Euclidean metric, depicted in Figure \ref{fig:ghngh} (we could equally well consider $\RR \times M \subseteq \RR^{1,2}$, there is nothing special about the negative definite case). We consider two different time frame choice, namely (A) $(E_1,E_2)$ for which the removed positive halfline $E_1+E_2$ is timelike, and (B) $(-E_1,E_2)$ for which it is ``spacelike'' (which we have not defined, meaning ``not causal''). For (A) this yields a non-globally hyperbolic spacetime because not all causal diamonds are compact. On other hand, for (B) all causal diamonds are compact because the removed line does not interfere, hence the resulting spacetime in the second case is globally hyperbolic.
 
 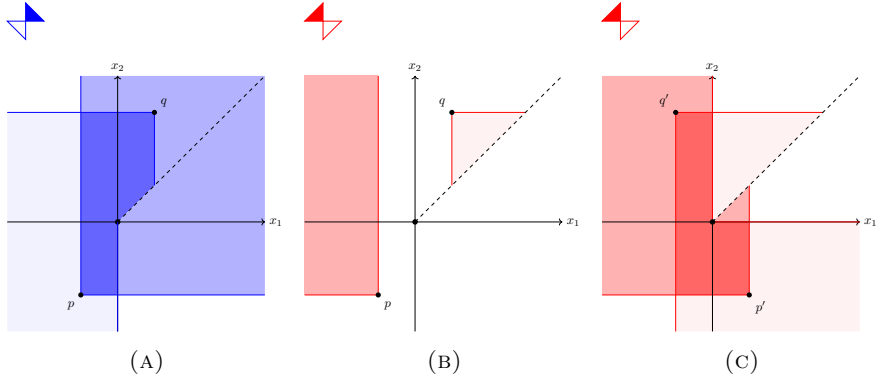
\begin{figure}[h]
  \begin{subfigure}[b]{0.3\textwidth}
   \centering
   \resizebox{\linewidth}{!}{
     \begin{tikzpicture}
  \coordinate (o) at (0,0);
  \coordinate (e) at (4,4);
  \coordinate (p) at (-1,-2);
  \coordinate (q) at (1,3);
  
   \fill[blue!60] (p) -- (-1,3) -- (q) -- (1,1) -- (o) -- (0,-2) -- cycle;
   \fill[blue!30] (-1,4) -- (-1,3) -- (q) -- (1,1) -- (o) -- (0,-2) -- (4,-2) -- (4,4) -- cycle;
   \fill[blue!5] (-3,3) -- (-1,3) -- (p) -- (0,-2) -- (0,-3) -- (-3,-3) -- cycle;
  \draw[-,draw=blue] (-3,3) -- (q) -- (1,1);
  \draw[-,draw=blue] (-1,4) -- (p) -- (4,-2);
 
  \draw[-,draw=blue] (-3,5.5) -- (-2.5,5) -- (-2.5,5.5) -- cycle;
  \filldraw[fill=blue,draw=blue] (-2.5,5.5) -- (-2.5,6) -- (-2,5.5) -- cycle;
  
    \draw[->] (-3,0) -- (4,0) node[right] {$x_1$};
  \draw[->] (0,-3) -- (0,4) node[above] {$x_2$};
   \draw[dashed] (o) node[circle,fill,inner sep=1.5pt]{}
   (o) -- (e);
     \draw[-,draw=blue] (o) -- (0,-3);
   \node[dot,label=above right:$q$] at (q) {};
   \node[dot,label=below left:$p$] at (p) {};
 \end{tikzpicture}   
   }
   \caption{} 
  \end{subfigure}
  \begin{subfigure}[b]{0.3\textwidth}
  \centering
   \resizebox{\linewidth}{!}{
      \begin{tikzpicture}
  \coordinate (o) at (0,0);
  \coordinate (e) at (4,4);
  \coordinate (p) at (-1,-2);
  \coordinate (q) at (1,3);
  
  \fill[red!30] (-3,-2) -- (p) -- (-1,4) -- (-3,4) -- cycle;
  \fill[red!5] (1,1) -- (q) -- (3,3) -- cycle;
  \draw[-,draw=red] (1,1) -- (q) -- (3,3); 
  \draw[-,draw=red] (-3,-2) -- (p) -- (-1,4);
 
  \draw[-,draw=red] (-2,5.5) -- (-2.5,5) -- (-2.5,5.5) -- cycle; 
  \filldraw[fill=red,draw=red] (-2.5,6) -- (-3,5.5) -- (-2.5,5.5) -- cycle;
  
    \draw[->] (-3,0) -- (4,0) node[right] {$x_1$};
  \draw[->] (0,-3) -- (0,4) node[above] {$x_2$};
   \draw[dashed] (o) node[circle,fill,inner sep=1.5pt]{}
   (o) -- (e);
   \node[dot,label=above left:$q$] at (q) {};
   \node[dot,label=below right:$p$] at (p) {};
 \end{tikzpicture}
   }
   \caption{} 
  \end{subfigure}
  \begin{subfigure}[b]{0.3\textwidth}
  \centering
  \resizebox{\linewidth}{!}{
     \begin{tikzpicture}
  \coordinate (o) at (0,0);
  \coordinate (e) at (4,4);
  \coordinate (p') at (1,-2);
  \coordinate (q') at (-1,3);
  
  \fill[red!30] (-3,-2) -- (p') -- (1,1) -- (o) -- (0,4) -- (-3,4) -- cycle;
  \fill[red!5] (-1,-3) -- (q') -- (3,3) -- (o) -- (4,0) -- (4,-3) -- cycle;
  \fill[red!60] (q') rectangle (0,-2); 
  \fill[red!60] (o) rectangle (p');
  \draw[-,draw=red] (-1,-3) -- (q') -- (3,3);
  \draw[-,draw=red] (-3,-2) -- (p') -- (1,1);
 
  \draw[-,draw=red] (-2,5.5) -- (-2.5,5) -- (-2.5,5.5) -- cycle;
  \filldraw[fill=red,draw=red] (-2.5,6) -- (-3,5.5) -- (-2.5,5.5) -- cycle;
  
    \draw[->] (-3,0) -- (4,0) node[right] {$x_1$};
  \draw[->] (0,-3) -- (0,4) node[above] {$x_2$};
     \draw[-,draw=red] (o) -- (4,0);
     \draw[-,draw=red] (o) -- (0,4);
        \draw[dashed] (o) node[circle,fill,inner sep=1.5pt]{}
   (o) -- (e);
   \node[dot,label=above left:$q'$] at (q') {};
   \node[dot,label=below right:$p'$] at (p') {};
 \end{tikzpicture}
 }
    \caption{} 
  \end{subfigure}
  \caption{\small We compare two different time frames on the same manifold $M = \RR^{0,2} \setminus (\RR_+\cdot(E_1+E_2))$, described in Example \ref{ex:ghngh}. In (A) the choice is $(E_1,E_2)$, in (B) and (C) we use $(-E_1,E_2)$. 
  The future is darker shaded than the past, the intersection (the causal diamond) is darkest.
  We see that spacetime (A) is not globally hyperbolic, because the  the causal diamond $J^+(p) \cap J^-(q)$ is not compact. In (B) this causal diamond is empty. Causal diamonds for other points $p'$ and $q'$ are also compact, as shown exemplary in (C), because the removed line is ``spacelike'' in this setting. This different behavior justifies the necessity of frames in our Definition \ref{def:sRspacetime} of $(n-\nu,\nu)$-spacetimes.}
  \label{fig:ghngh}
 \end{figure}
\end{ex}

\subsection{Global hyperbolicity iterated}\label{ssec:sRGH}

 In Theorem~\ref{prop:equivLR} (i) we have related the Heine--Borel property of a Riemannian manifold $(\Sigma,\sigma)$ \emph{directly} to the compactness of causal diamonds of the corresponding Lorentzian product spacetime $(\RR \times \Sigma, -dt^2 \oplus \sigma)$. Can we prove a similar connection for the related notions of global hyperbolicity for different $\nu$?
 
 Since we have already treated the $\nu=0$ case we can restrict to $0 < \nu \leq n$. We first attempt a naive approach. By that we mean that to go from a spacetime structure on $(M,g)$ given by vector fields $X_1,\ldots,X_\nu$ we define one on the product
\begin{align}\label{Mprime}
 (M',g') = (\RR \times M, -dt^2 \oplus g)
\end{align}
 by using the additional vector field $X_{\nu+1} =\partial_t$. We call this $(n-\nu,\nu+1)$-spacetime structure on $(M',g')$ the \emph{orthogonally extended $(n-\nu,\nu+1)$-spacetime}. We also consider the product
 \begin{align}\label{Mdprime}
  (M'',g'') = (\RR \times M, dt^2 \oplus g)
 \end{align}
 equipped with the same spacetime structure but viewed as a $(n+1-\nu,\nu)$-spacetime.
 
 We show that the answer to the above question is no, i.e., that the notions of global hyperbolicity for different $\nu$ are in general unrelated. For instance, we demonstrate that global hyperbolicity of $(M,g)$ does not imply global hyperbolicity of $(M',g')$. On the other hand, we succeed in showing that the top-down implication for $(M',g')$ to $(M,g)$ still holds. In this section we prove Theorem \ref{thm:equivsR} and provide a counterexample for the missing bottom-up implication. The reason why going from $(M,g)$ to $(M',g')$ does not work is, in short, because the boundary $\partial C'$ of the corresponding cone structure includes all of $C$ and is ``too wide''. In other words, $M'$ leaves more room than in $M$ because $g'(v,v) \leq 0$ does not imply $g(v,v) \leq 0$ (this was not an issue for positive definite $g$). The metric perspective and the fact that we crucially exploited the direct (and not just conformal) relation between $d_\sigma$ and $\hat d_t$ in Section \ref{sec:LRequiv} also offers an explanation for why our proof from Theorem~\ref{prop:equivLR} cannot be extended.
 
\medskip
 Before jumping right into the analysis of global hyperbolicity in this framework it is insightful and necessary to consider (stable) causality first. Indeed, we see that problems already occur at this level.
 
 \begin{lem}\label{lem:Cs}
  Let $(M,g)$ be a $(n-\nu,\nu)$-spacetime with $0\leq \nu \leq n$. Then $(M',g')$ and $(M'',g'')$ as defined in \eqref{Mprime} and \eqref{Mdprime} are $(n-\nu,\nu+1)$- and $(n+1-\nu,\nu)$-spacetimes, respectively. The corresponding cone structures satisfy
  \begin{align}\label{Cs}
   C'' \preccurlyeq \RR \times C \quad \text{and} \quad \{0\} \times C \preccurlyeq C' \cap C''.
  \end{align}
  If $\nu = n$ then furthermore
  \begin{align}\label{Cs2}
    C' = ([0,\infty) \times C) \cup ((0,\infty) \times TM_0). 
  \end{align}
 \end{lem}

 \begin{proof}
  Since all cone structures are continuous by Theorem~\ref{thm:propercone} it is sufficient to prove the set-theoretic implications pointwise. Since the map $A \mapsto \overline{\operatorname{co}}~A$ is nondecreasing it is furthermore sufficient to prove the statement for the closed precausal cones $\mathcal{C}$ used in the proof of Theorem~\ref{thm:propercone}. We first prove the inclusions of \eqref{Cs}. Suppose the time frame orientation of $M$ is given by the vector fields $X_1,\ldots,X_\nu$, and $X_{\nu+1}=\partial_t$. We denote the pulled back version onto $\RR \times M$ along $\pi_M$ in the same way.
  
  If $v = (v_0,v_M) \in \mathcal{C}'' \subseteq TM'' = \RR \times TM$ then
  \begin{align*}
   g(v_M,v_M) & = g''(v,v) - dt^2(v_0,v_0) \leq 0, \\
   g(v_M,X_i) & = g''(v,X_i) - \underbrace{g''(v_0,X_i)}_{=0} \leq 0,
  \end{align*}
  with a strict inequality for at least one $i=1,\ldots,\nu$ in the last line. Hence $v_M \in \mathcal{C}' \subseteq TM \setminus (TM)_0$ irrespective of the value of $v_0$. In other words, $\mathcal{C}'' \preccurlyeq \RR \times \mathcal{C}$.
  
  If, on the other hand, $v_M \in \mathcal{C} \subseteq TM \setminus (TM)_0$ then $v = (0,v_M) \neq 0$ satisfies due to the orthogonal product structure
  \begin{align*}
   g'(v,v) \leq g''(v,v) &= g(v_M,v_M) \leq 0, \\
   g'(v,X_i) = g''(v,X_i) &= g(v_M,X_i) \leq 0, \qquad i=1,\ldots,\nu.
  \end{align*}
  Moreover, $g'(v,\partial_t) = 0$. Thus $v \in \mathcal{C}' \cap \mathcal{C}''$.
  
  It remains to prove \eqref{Cs2} for $\nu = n$. Consider once more $v=(v_0,v_M) \in \mathcal{C}'=C'$. Due to the assumption $g'(v,\partial_t)\leq 0$ we must have $v_0 \geq 0$. In addition, for all $i=1,\ldots,n$, we have
  \begin{align}\label{vM1}
   g(v_M,X_i) = g'(v,X_i) \leq 0.
  \end{align}
  Since $g$ is negative definite by assumption, we have for all $v_M$
  \begin{align}\label{vM2}
   g(v_M,v_M) \leq 0,
  \end{align}
  with equality if and only if $v_M=0$. If $v_M \neq 0$ then \eqref{vM1}--\eqref{vM2} show that $v_M \in \mathcal{C}=C$. If, on the other hand, $v_M=0$ then by Lemma \ref{lem:onneg} we must have $g'(v,\partial_t)< 0$, i.e., $v_0 > 0$. We have thus shown the inclusion $\preccurlyeq$ in \eqref{Cs2}. The other inclusion $\succcurlyeq$ follows in the same way from the above inequalities and cases.
 \end{proof}
 
 \begin{rem}[$C'$ and $C''$]
 Note that $C' \not\preccurlyeq C''$ because $\partial_t \in C' \setminus C''$. On the other hand, for a unital vector field $X \in C$ we have $Y = X -\partial_t \in C'' \setminus C'$ because $g'(Y,\partial_t) = 1$.
 \end{rem}
 
 The implications of \eqref{Cs} allow us to efficiently relate the causal structures of $M$, $M'$, and $M''$, which ultimately cumulates to the prove of Theorem \ref{thm:equivsR}.

 \begin{prop}\label{prop:causal}
  Let $(M,g)$ be a $(n-\nu,\nu)$-spacetime with $0\leq\nu\leq n$. Then for the spacetimes $(M',g')$ and $(M'',g'')$ as defined in \eqref{Mprime} and \eqref{Mdprime} we have
  \begin{align*}
   & (M',g') \text{ is causal} 
   \Longleftrightarrow (M,g) \text{ is causal} 
   \Longleftrightarrow (M'',g'') \text{ is causal}.
  \end{align*}
 \end{prop}

 \begin{proof}
  ($M' \Longrightarrow M \Longleftarrow M''$) By \eqref{Cs}, if there is a closed causal curve $\gamma_M$ in $M$, then $\gamma = (0,\gamma_M)$ is a closed causal curve in both $M'$ and $M''$. Thus causality of $M'$ (or $M''$) enforces causality of $M$.
  
  ($M \Longrightarrow M''$) If $\gamma$ is a closed $g''$-causal curve in $M''$, then by \eqref{Cs} the projection $\gamma_M = \pi_M \circ \gamma$ is a closed causal curve in $M$, a contradiction.
  
  ($M' \Longleftarrow M$) Suppose $\gamma$ is a closed $g'$-causal curve in $M'$. Without loss generality we assume that it is future directed, thus either (i) $\dot\gamma_0 = 0$ almost everywhere or (ii) $\dot\gamma_0 \neq 0$ on a set of nonzero Lebesgue measure. In the case of (i), since $\gamma_0$ is absolutely continuous and the fundamental theorem of calculus applies, $\gamma_0$ must be constant and thus the projection $\gamma_M$ is a closed $g$-causal curve in $M$, a contradiction. In case of (ii) there must be a set of nonzero Lebesgue measure where $\dot\gamma_0 >0$ and thus, due to closedness of the curve in the $t$ variable and the fundamental theorem of calculus, there must also be a set of nonzero Lebesgue measure where $\dot\gamma_0<0$, a contradiction to $\gamma$ being future directed causal. Thus $\gamma$ does not exist.
 \end{proof}
 
 Note that Proposition \ref{prop:causal} is already weaker than in the Riemannian vs.\ Lorentzian case ($\nu=0$) where $M'$ is (even stably!) causal completely independent of $M$. From the proof it is also clear which implication is the most fragile one. Indeed, we lose it when stepping up on the causal ladder.
 
 \begin{prop}\label{prop:Kcausal}
  Let $(M,g)$ be a $(n-\nu,\nu)$-spacetime with $0\leq \nu \leq n$, and let $(M',g')$ and $(M'',g'')$ be as defined in \eqref{Mprime} and \eqref{Mdprime}, respectively. Then
  \begin{align*}
   & (M',g') \text{ is stably causal} \\
   \Longrightarrow \, & (M,g) \text{ is stably causal} \\
   \Longleftrightarrow \, & (M'',g'') \text{ is stably causal}.
  \end{align*}
  To be precise, for each valid implication, the same (up to natural projection and embedding) time functions can be used.
  
  Moreover, if $\nu = n$ and $\tau_M$ is a time function for $M$ then the function $\tau = \tau_M \circ \pi_M + t$ is a time function for $M'$, which also shows that
  \begin{align*}
   (M',g') \text{ is stably causal}
   \Longleftarrow (M,g) \text{ is stably causal}.
  \end{align*}
 \end{prop}
 
 \begin{proof}
  ($M'\Longrightarrow M \Longleftarrow M''$) Since the cone structures satisfies $\{0\} \times C \preccurlyeq C' \cap C''$ by \eqref{Cs} the causal curves in $M$ are fully captured by those in $C'$ (and in $C''$). Thus a time function $\tau$ of $M'$ (or of $M''$) induces a time function $\tau_M(p_M) := \tau(0,p_M)$ on $M$.
  
  ($M \Longrightarrow M''$) Suppose $\tau_M$ is a time function of $M$. Since by \eqref{Cs} $\pi_M(C'') \preccurlyeq C$ we can simply pullback $\tau_M$ to $M''$ to obtain a time function $\tau = \tau_M \circ \pi_M$. 

  ($M \Longrightarrow M'$ if $\nu=n$) Suppose $\tau_M$ is a time function for $M$ and $t$ the coordinate of the additional dimension. It follows immediately from \eqref{Cs2} that $\tau_M \circ \pi_M + t$ is a time function for $M'$.
 \end{proof}

 Recall that the reason why the second implication ($M \Longrightarrow M''$) in Proposition \ref{prop:Kcausal} holds is that compared to $M$ the projections of the cones in $M''$ to $M$ are narrower (due to more positive directions) or remain the same (if the new variable remains fixed). On the other hand, the projections of the cones in $M'$ become strictly wider than the cones on $M$ (due to more negative tangent directions) even if the new variable is fixed (unless $\nu =n$, in which case they remain the same). In other words, if $\gamma$ is a future directed $g'$-causal curve then $\pi_M \circ \gamma$ need not be future directed $g$-causal. This also explains why the implication ($M \Longrightarrow M'$) only works for $\nu=n$.
 
 Moreover, unlike from the Riemannian to the Lorentzian case the coordinate $t$ is not a time function for $M'$ if $X_{\nu+1} = \partial_t$. We show that if one mildly perturbs the vector field $\partial_t$ one can nonetheless ensure that $t$ is a time function and thus safe stable causality for $M'$ while leaving the first $\nu$ time frame orienting vector fields intact. As in Section~\ref{ssec:Lprod} (stable) causality of $M$ is then also not required.
 
\begin{prop}\label{prop1stepback}
 Let $(M,g)$ be a $(n-\nu,\nu)$-spacetime with $0< \nu\leq n$ and time frame orientation defining vector fields $X_1,\ldots,X_\nu$. Then for every $\varepsilon >0$ the semi-Riemannian manifold $(M',g')=(\RR\times M, -dt^2 \oplus g)$ with vector fields $X_1,\ldots,X_\nu$ and $X_{\nu+1} = \partial_t - \varepsilon (X_1 + \ldots + X_\nu)$ is a $(n-\nu,\nu+1)$-spacetime with temporal function $t$. In particular, $(M',g')$ is stably causal.
\end{prop}

\begin{proof}
 Since the vector fields $X_1,\ldots,X_\nu,X_{\nu+1}$ are continuous, linearly independent, and satisfy
 \begin{align*}
  g'(X_i,X_i) &= g(X_i,X_i) < 0, \qquad i =1,\ldots,\nu, \\
  g'(X_{\nu+1},X_{\nu+1}) &= - 1 + \varepsilon^2 \underbrace{g(X_1+\ldots + X_\nu, X_1+\ldots + X_\nu)}_{< 0} < 0,
 \end{align*}
 they define a time frame orientation for $(M',g')$ that is based on the time frame orientation of the submanifolds $\{p_0\}\times M$ (but not extended orthogonally!).
 
 We  show that $t$ is a temporal function. To this end write any $v \in C' \subseteq T_pM'$ as direct sum $v = v_0 + v_M \in \spann (\partial_t(p)) \oplus \spann (\partial_t(p))^\perp$. Future directedness in the first $\nu$ directions and orthogonality $g' = -dt^2 \oplus g$ implies that
 \begin{align}\label{vMcond}
  g(v_M,X_i(p_M)) = g'(v,X_i(p)) \leq 0, \qquad i=1,\ldots,\nu. 
 \end{align}
 In particular, $g(v_M,X_1+\ldots + X_\nu) \leq 0$.
 Future directedness of $v$ in the $(\nu+1)$-th direction means that
 \begin{align}\label{dtineq}
  0 \geq g'(v,X_{\nu+1}) &= g'(v,\partial_t) - \varepsilon g'(v,X_1+\ldots + X_\nu) \\
                  &= - dt(v) - \varepsilon g(v_M,X_1+\ldots + X_\nu). \nonumber
 \end{align}
 Note that, in addition to \eqref{vMcond} and \eqref{dtineq}, by Lemma~\ref{lem:onneg} there is at least one $j \in \{1,\ldots,\nu+1\}$ such that
 \[
   g'(v,X_j(p)) < 0.
 \]
 If this is the case for $j = \nu+1$ then \eqref{dtineq} is a strict inequality and hence
 \[
  dt(v) > - \varepsilon g(v_M,X_1+\ldots + X_\nu) \geq 0.
 \]
 Otherwise, \eqref{dtineq} is an equality but there must be $j \in \{1,\ldots,\nu \}$ for which \eqref{vMcond} is a strict inequality. Hence
 \[
  dt(v) = - \varepsilon g(v_M,X_1+\ldots + X_\nu) \geq - \varepsilon g(v_M,X_j) > 0.
 \]
 Either way, for all future directed causal vectors $v \in T_pM$ we have $dt(v) >0$. Hence $t$ is a temporal function of $(M',C')$.
\end{proof}
 
 Having established the equivalence of causality in Proposition~\ref{prop:causal} and having shown how time functions carry over from one spacetime to another in Propositions~\ref{prop:Kcausal} and \ref{prop1stepback}, we finally turn to Cauchy time functions and global hyperbolicity and prove the remaining implications of Theorem~\ref{thm:equivsR}.
 
 \medskip
 We start with the easier situation $M''$. Adding positive tangent directions, meaning $+dt^2$, is not a problem as no additional time frame orientation defining vector fields are needed and projections and pullbacks behave as desired.

\begin{prop}\label{prop4}
 Let $(M,g)$ be a $(n-\nu,\nu)$-spacetime with time frame orientation defining vector fields $X_1,\ldots,X_\nu$, and $(M'',g'') = (\RR \times M, dt^2 \oplus g)$ the $(n-\nu+1,\nu)$-spacetime via the same vector fields. Then
 \[ (M,g) \text{ globally hyperbolic } \Longleftrightarrow (M'',g'') \text{ globally hyperbolic}.\]
\end{prop}

\begin{proof}
 ($\Longrightarrow$) If $(M,g)$ is globally hyperbolic then, by Theorem~\ref{thm:ghtime}, there exists a Cauchy time function $\tau_M \colon M \to \RR$. In the proof of Proposition~\ref{prop:Kcausal} we have seen that $\tau = \tau_M \circ \pi_M$ is a time function for $M''$. 

 Suppose $\tau$ is not Cauchy. Then there is an inextendible causal curve $\alpha$ such that $\sup(\tau\circ\alpha) < \infty$. Since $\RR$ is complete we can always extend in the first component $\alpha_0$. Moreover, since $\alpha_M = \pi_M \circ \alpha$ is a causal curve in $M$ and $(M,g)$ is globally hyperbolic we have that $\alpha_M$ must be extendible in $M$, and hence $\alpha$ must be extendible in $M''$, a contradiction. Hence $\tau$ is Cauchy in $(M'',g'')$ which is therefore globally hyperbolic by Theorem~\ref{thm:ghtime}.
 
 ($\Longleftarrow$) Suppose $(M'',g'')$ is globally hyperbolic and $\tau$ is a Cauchy time function for $M''$. By Proposition \ref{prop:Kcausal} the pulled backed $\tau_M(p_M) = \tau(0,p_M)$ is a time function for $M$. If $\tau_M$ is not Cauchy then there is an inextendible future directed causal curve $\alpha_M \colon \RR \to M$ such that $\sup(\tau_M \circ \alpha_M)<\infty$. For the pulled back curve $\alpha=(0,\alpha_M)$ we would then have $\sup(\tau\circ\alpha)<\infty$, meaning that it is extendible in $M''=\RR\times M$ because of Cauchyness of $\tau$. Since the first component of $\alpha$ is constant this implies that $\alpha_M$ is extendible in $M$. Hence $\tau_M$ is a Cauchy time function.
\end{proof}

\begin{prop}\label{prop1step}
 Let $(M,g)$ be a $(n-\nu,\nu)$-spacetime with $0 < \nu \leq n$ and corresponding proper cone structure $(M,C)$. Then for the orthogonally extended $(n-\nu,\nu+1)$-spacetime $(M',g') = (\RR \times M, - dt^2 \oplus g)$ we have that 
 \[ (M',g') \text{ globally hyperbolic } \Longrightarrow (M,g)  \text{ globally hyperbolic}.
 \]
 If $\nu = n$ then also the reverse implication ($\Longleftarrow$) holds.
\end{prop}

\begin{proof}
 By Lemma~\ref{lem:Cs} $(M',g')$ is a spacetime with $\{ 0 \} \times C \preccurlyeq C'$. That global hyperbolicity descends from $M'$ to $M$ thus follows as in the proof of Proposition \ref{prop4} for $M''$.
 
 It remains to prove the reverse implication if $\nu = n$. If $\tau_M$ is a Cauchy time function in $M$ then we already know from Proposition \ref{prop:Kcausal} that $\tau = \tau_M \circ \pi_M + t$ is a time function for $M'$. Suppose $\tau$ is not Cauchy in $M'$. Then there is a future directed future inextendible $g'$-causal curve $\alpha$ such that $\sup(\tau \circ \alpha) < \infty$. This implies that both $\tau_M \circ \alpha_M$ and $t \circ \alpha$ are bounded from above, but since $M$ is globally hyperbolic (and $\dot\alpha_M$ is $g$-causal or $=0$ by \eqref{Cs2}) and $\mathbb{R}$ is complete, we can extend both $\alpha_M$ and $\alpha_0$ to the future, and hence $\alpha$ as well, a contradiction.
\end{proof}

 We conclude with an example that demonstrates that, in general, the implication of Proposition \ref{prop1step} cannot be reversed if $\nu<n$. This should come as no surprise because in Lemma~\ref{lem:Cs} we have observed that the cones in $M'$ include the base $M$ if we choose $X_{\nu+1}=\partial_t$. However, we argue that even when turning them into true cones, for instance, by using $X_{\nu+1} = \partial_t - \varepsilon (X_1 + \ldots + X_\nu)$ for $\varepsilon >0$ as in Proposition~\ref{prop1stepback} for enforcing stable causality on $M'$, global hyperbolicity of $M$ does not, in general, imply global hyperbolicity of $M'$. This is in stark contrast to the Riemannian and Lorentzian equivalence obtained in Theorem~\ref{prop:equivLR}.

\begin{ex}[Global hyperbolicity of $M'$ not inherited from $M$ if $\nu < n$]\label{ex:notgh}
 Let $M = \RR^{1,1} \setminus J^+(0) \subseteq \RR^{1,1}$ be the a submanifold of $2$-dimensional Minkowski space with metric tensor $\eta = -dx^2 + dy^2$ (but same holds if $\dim M >2$). It is a spacetime with time orientation defining vector field $X_1 = \partial_x$. Since $M$ is causal and all causal diamonds are compact it is easy to see that $M$ is globally hyperbolic. In particular, for the points $p_M = (-1,1)$ and $q_M=(1,-2)$ we have $J^+(p_M) \cap J^-(q_M) = \emptyset$. See also Figure \ref{fig:M}.
 
 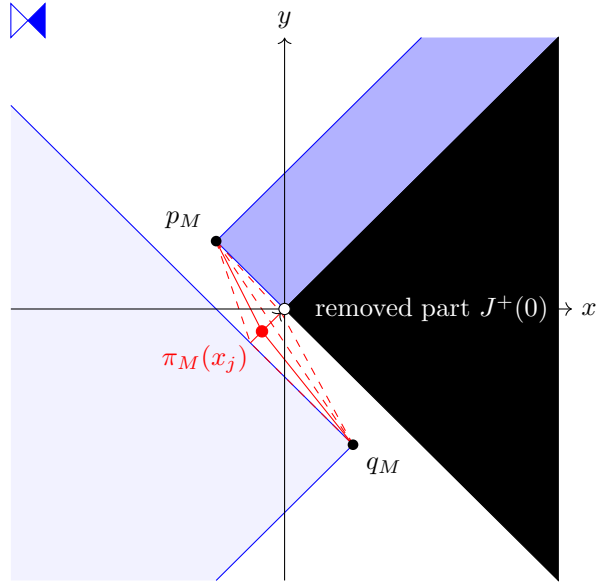
\begin{figure}[h]
    \begin{tikzpicture}[scale=0.9]
  \coordinate (o) at (0,0);
  \coordinate (e) at (4,4);
  \coordinate (w) at (-4,-4);
  \coordinate (pM) at (-1,1);
  \coordinate (qM) at (1,-2);
  
  \draw[-,draw=blue] (-4,4.5) -- (-3.75,4.25) -- (-4,4) -- cycle;
  \filldraw[fill=blue,draw=blue] (-3.75,4.25) -- (-3.5,4.5) -- (-3.5,4) -- cycle;
  
  \fill[blue!30] (2,4) -- (pM) -- (o) -- (e) -- cycle;
  \draw[-,draw=blue] (2,4) -- (pM) -- (o);
  
  \fill[blue!5] (-4,3) -- (qM) -- (-1,-4) -- (w) -- cycle;
  \draw[-,draw=blue] (-4,3) -- (qM) -- (-1,-4);

    \draw[dashed,draw=red] (pM) -- (-0.5,-0.5) -- (qM);
    \draw[-,draw=red] (pM) -- (-0.33,-0.33) node[circle,fill=red,draw=red,inner sep=1.5 pt,label={[text=red]below left:$\pi_M(x_j)$}] {} (-0.33,-0.33) -- (qM);
    
    \draw[dashed,draw=red] (pM) -- (-0.2,-0.2) -- (qM);
    \draw[dashed,draw=red] (pM) -- (-0.05,-0.05) -- (qM);
    \draw[->,draw=red] (-0.5,-0.5) -- (-0.05,-0.05);
  
    \draw[->] (-4,0) -- (4.2,0) node[right] {$x$};
  \draw[->] (0,-4) -- (0,4) node[above] {$y$};
    \filldraw[-,fill=black] (e) -- (o) -- (4,-4) -- cycle;
    \draw[-] (o) node[circle,fill=white,draw=black,inner sep=1.5pt,label={[text=white]right:$\;\;$removed part $J^+(0)$}]{}
   (o) -- (e);
   \node[dot,label=below right:$q_M$] at (qM) {};
   \node[text=blue,align=left] at (-2,-2) {$J^-(q_M)$};
   \node[dot,label=above left:$p_M$] at (pM) {};
   \node[text=blue,align=right] at (1,2) {$J^+(p_M)$};
   \end{tikzpicture}
   \caption{\small The figure depicts the projection of $M'$ onto $M$ as used in Example \ref{ex:notgh}. The spacetime $M$ is globally hyperbolic, however $M'=\RR\times M$ is not. If one considers the points $p$ and $q$ one can see that in $M$ the causal diamond $J^+(p_M) \cap J^-(q_M) = \emptyset$ is compact. On the other hand, in $M'$ the causal diamond $J^+(p) \cap J^-(q)$ is not compact, because the red paths $\alpha_j$ and $\beta_j$ connecting $p$ and $q$ with $x_j$, respectively, are causal, however the sequence $(x_j)$ does not have a limit because $0 \not\in M'$.}
   \label{fig:M}
 \end{figure}

 Let us consider $(M',g') = (\RR \times M, -dt^2 \oplus \eta)$ with additional vector field $X_2 = \partial_t$. We will show that $(M',g')$ is not globally hyperbolic: Consider $p = (-2, p_M)$ and $q=(2,q_M)$. Then $J^+(p) \cap J^-(q) \neq \emptyset$ but the causal diamond not closed: Consider the sequence of points $x_j = (0,-\frac{1}{j},-\frac{1}{j})$. Then the curves $\alpha_j \colon [0,1] \to M'$,
 \[\alpha_j(s) = \left(2s-2,-1+s-\frac{s}{j},1-s- \frac{s}{j}\right),\]                                                                            
 are from $p$ to $x_j$. Since $\dot\alpha_j=(2,1-\frac{1}{j},-1-\frac{1}{j})$ one obtains for $j>1$
 \begin{align*}
  &g'(\dot\alpha_j,\dot\alpha_j) = -4 + \frac{4}{j} < 0, \\ &g'(\dot\alpha_j,X_1) = -1+\frac{1}{j}<0, \qquad g'(\dot\alpha_j,X_2) = -2 < 0,
 \end{align*}
 i.e., the curves are future directed (pre)causal in $M'$ (but not in $M$). Thus $x_j \in J^+(p)$ for $j>1$. Moreover, the curves $\beta_j \colon [0,1] \to M'$,
 \[\beta_j(s) = \left(2-2s,1-s-\frac{s}{j},-2+2s-\frac{s}{j}\right)\] are past directed (pre)causal in $M'$ from $q$ to $x_j$ for all $j\geq 1$, because $\dot\beta_j = (-2,-1-\frac{1}{j},2-\frac{1}{j})$ and thus
 \begin{align*}
  & g'(\dot\beta_j,\dot\beta_j) = -1 - \frac{6}{j} < 0, \\
  & g'(\dot\beta_j,X_1) = 1 + \frac{1}{j} >0, \qquad g'(\dot\beta_j,X_2) = 2>0.
 \end{align*}
 Thus $x_j \in J^-(q)$ for all $j$. In particular, $x_j \in J^+(p) \cap J^-(q)$ for all $j \geq 2$ and the causal diamond is nonempty. The limiting point of this sequence, $x = \lim_{j\to \infty}x_j = (0,0,0)$, however, is not in $M'$. Thus $J^+(p) \cap J^-(q)$ is not closed and hence also not compact. Therefore, by definition, $(M',g')$ is not a globally hyperbolic spacetime.

 One could have the hope that by replacing $X_2$ by $X_2 -\varepsilon X_1$ for appropriately chosen $\varepsilon$ one could still carry over global hyperbolicity from $M$ to $M'$. This is generally not the case. One can still find examples of the type studied above as $\partial_t$ always remains to be a future directed causal vector field for any choice of $\varepsilon$ in the assumptions of Proposition~\ref{prop1stepback} (even if it is negative) as $g' = -dt^2 \oplus g$ consists of orthogonal components and one can go infinitely fast in the $t$-direction to compensate for any $\varepsilon$. Note that such examples can only be constructed for $0<\nu < n$ and the problem may very well be that the Wick-rotated metric $dx^2+dy^2$ is \emph{not} a complete Riemannian metric on $M$.
\end{ex}

 All in all, we have shown Theorem~\ref{thm:equivsR}. The implications are established in Propositions~\ref{prop4} and \ref{prop1step}, and Example~\ref{ex:notgh} shows that the missing implication is indeed false.


\section*{Acknowledgments}
 
The author would like to thank Jim Isenberg for encouraging her to explore the semi-Riemannian setting and Steffen Sagave for discussing some references in algebraic topology. This project was supported by the Dutch Research Council (NWO), Project number VI.Veni.192.208.
%

\section*{Statements and declarations}

 This work was supported by the Dutch Research Council (NWO), Project number VI.Veni.192.208. The author states that there is no conflict of interest. No data was generated or used in this manuscript.


\bibliographystyle{abbrv}
\bibliography{refs.bib}

\begin{bibdiv}
\begin{biblist}

\bib{Ada:spherevfs}{article}{
      author={Adams, J.~F.},
       title={Vector fields on spheres},
        date={1962},
        ISSN={0003-486X},
     journal={Ann. of Math. (2)},
      volume={75},
       pages={603\ndash 632},
         url={https://doi.org/10.2307/1970213},
      review={\MR{139178}},
}

\bib{Al1}{article}{
      author={Allen, Brian},
       title={Null distance and {G}romov-{H}ausdorff convergence of warped
  product spacetimes},
        date={2023},
        ISSN={0001-7701,1572-9532},
     journal={Gen. Relativity Gravitation},
      volume={55},
      number={10},
       pages={Paper No. 118, 34},
         url={https://doi.org/10.1007/s10714-023-03167-8},
}

\bib{AB}{article}{
      author={Allen, Brian},
      author={Burtscher, Annegret},
       title={Properties of the null distance and spacetime convergence},
        date={2022},
        ISSN={1073-7928},
     journal={Int. Math. Res. Not. IMRN},
      number={10},
       pages={7729\ndash 7808},
         url={https://doi.org/10.1093/imrn/rnaa311},
      review={\MR{4418719}},
}

\bib{AtDu}{article}{
      author={Atiyah, M.~F.},
      author={Dupont, J.~L.},
       title={Vector fields with finite singularities},
        date={1972},
        ISSN={0001-5962,1871-2509},
     journal={Acta Math.},
      volume={128},
       pages={1\ndash 40},
         url={https://doi.org/10.1007/BF02392157},
      review={\MR{451256}},
}

\bib{At}{book}{
      author={Atiyah, Michael~F.},
       title={Vector fields on manifolds},
      series={Arbeitsgemeinschaft f\"{u}r Forschung des Landes
  Nordrhein-Westfalen},
   publisher={Westdeutscher Verlag, Cologne},
        date={1970},
      volume={Heft 200},
      review={\MR{263102}},
}

\bib{AuCe:setvalued}{book}{
      author={Aubin, Jean-Pierre},
      author={Cellina, Arrigo},
       title={Differential inclusions},
      series={Grundlehren der mathematischen Wissenschaften [Fundamental
  Principles of Mathematical Sciences]},
   publisher={Springer-Verlag, Berlin},
        date={1984},
      volume={264},
        ISBN={3-540-13105-1},
         url={https://doi.org/10.1007/978-3-642-69512-4},
        note={Set-valued maps and viability theory},
      review={\MR{755330}},
}

\bib{BEE}{book}{
      author={Beem, John~K.},
      author={Ehrlich, Paul~E.},
      author={Easley, Kevin~L.},
       title={Global {L}orentzian geometry},
     edition={Second},
      series={Monographs and Textbooks in Pure and Applied Mathematics},
   publisher={Marcel Dekker, Inc., New York},
        date={1996},
      volume={202},
        ISBN={0-8247-9324-2},
}

\bib{BL:3mf}{article}{
      author={Benedetti, Riccardo},
      author={Lisca, Paolo},
       title={Framing 3-manifolds with bare hands},
        date={2018},
        ISSN={0013-8584,2309-4672},
     journal={Enseign. Math.},
      volume={64},
      number={3-4},
       pages={395\ndash 413},
         url={https://doi.org/10.4171/LEM/64-3/4-9},
      review={\MR{3987151}},
}

\bib{BeSa:orthsplit}{article}{
      author={Bernal, Antonio~N.},
      author={S\'{a}nchez, Miguel},
       title={Smoothness of time functions and the metric splitting of globally
  hyperbolic spacetimes},
        date={2005},
        ISSN={0010-3616,1432-0916},
     journal={Comm. Math. Phys.},
      volume={257},
      number={1},
       pages={43\ndash 50},
         url={https://doi.org/10.1007/s00220-005-1346-1},
      review={\MR{2163568}},
}

\bib{BeSu:time}{article}{
      author={Bernard, Patrick},
      author={Suhr, Stefan},
       title={Lyapounov functions of closed cone fields: from {C}onley theory
  to time functions},
        date={2018},
        ISSN={0010-3616,1432-0916},
     journal={Comm. Math. Phys.},
      volume={359},
      number={2},
       pages={467\ndash 498},
         url={https://doi.org/10.1007/s00220-018-3127-7},
      review={\MR{3783554}},
}

\bib{BeSu:gh}{article}{
      author={Bernard, Patrick},
      author={Suhr, Stefan},
       title={Cauchy and uniform temporal functions of globally hyperbolic cone
  fields},
        date={2020},
        ISSN={0002-9939,1088-6826},
     journal={Proc. Amer. Math. Soc.},
      volume={148},
      number={11},
       pages={4951\ndash 4966},
         url={https://doi.org/10.1090/proc/15106},
      review={\MR{4143406}},
}

\bib{BDS:cobord}{article}{
      author={B\"{o}kstedt, Marcel},
      author={Dupont, Johan},
      author={Svane, Anne~Marie},
       title={Cobordism obstructions to independent vector fields},
        date={2015},
        ISSN={0033-5606,1464-3847},
     journal={Q. J. Math.},
      volume={66},
      number={1},
       pages={13\ndash 61},
         url={https://doi.org/10.1093/qmath/hau011},
      review={\MR{3356279}},
}

\bib{BDGS}{article}{
      author={Borde, A.},
      author={Dowker, H.~F.},
      author={Garcia, R.~S.},
      author={Sorkin, R.~D.},
      author={Surya, S.},
       title={Causal continuity in degenerate spacetimes},
        date={1999},
        ISSN={0264-9381,1361-6382},
     journal={Classical Quantum Gravity},
      volume={16},
      number={11},
       pages={3457\ndash 3481},
         url={https://doi.org/10.1088/0264-9381/16/11/303},
      review={\MR{1729954}},
}

\bib{Br:Sn}{article}{
      author={Brouwer, L.E.J.},
       title={{\"U}ber {A}bbildung von {M}annigfaltigkeiten},
    language={ger},
        date={1912},
     journal={Mathematische Annalen},
      volume={71},
       pages={97\ndash 115},
         url={http://eudml.org/doc/158520},
}

\bib{BBI}{book}{
      author={Burago, Dmitri},
      author={Burago, Yuri},
      author={Ivanov, Sergei},
       title={A course in metric geometry},
      series={Graduate Studies in Mathematics},
   publisher={American Mathematical Society, Providence, RI},
        date={2001},
      volume={33},
        ISBN={0-8218-2129-6},
         url={https://doi.org/10.1090/gsm/033},
}

\bib{BGH2}{article}{
      author={Burtscher, Annegret},
      author={Garc\'{\i}a-Heveling, Leonardo},
       title={Global hyperbolicity through the eyes of the null distance},
        date={2024},
        ISSN={0010-3616,1432-0916},
     journal={Comm. Math. Phys.},
      volume={405},
      number={4},
       pages={Paper No. 90, 35},
         url={https://doi.org/10.1007/s00220-024-04936-5},
}

\bib{BGH1}{article}{
      author={Burtscher, Annegret},
      author={Garc\'ia-Heveling, Leonardo},
       title={Time functions on {L}orentzian length spaces},
        date={2025},
        ISSN={1424-0637,1424-0661},
     journal={Ann. Henri Poincar\'e},
      volume={26},
      number={5},
       pages={1533\ndash 1572},
         url={https://doi.org/10.1007/s00023-024-01461-y},
}

\bib{Bur}{article}{
      author={Burtscher, Annegret~Y.},
       title={Length structures on manifolds with continuous {R}iemannian
  metrics},
        date={2015},
        ISSN={1076-9803},
     journal={New York J. Math.},
      volume={21},
       pages={273\ndash 296},
         url={http://nyjm.albany.edu:8000/j/2015/21_273.html},
      review={\MR{3358543}},
}

\bib{Ca89}{article}{
      author={Carri\`ere, Yves},
       title={Autour de la conjecture de {L}. {M}arkus sur les vari\'{e}t\'{e}s
  affines},
        date={1989},
        ISSN={0020-9910,1432-1297},
     journal={Invent. Math.},
      volume={95},
      number={3},
       pages={615\ndash 628},
         url={https://doi.org/10.1007/BF01393894},
}

\bib{CrSt}{article}{
      author={Crabb, M.~C.},
      author={Steer, B.},
       title={Vector-bundle monomorphisms with finite singularities},
        date={1975},
        ISSN={0024-6115,1460-244X},
     journal={Proc. London Math. Soc. (3)},
      volume={30},
       pages={1\ndash 39},
         url={https://doi.org/10.1112/plms/s3-30.1.1},
      review={\MR{370614}},
}

\bib{Du:K}{article}{
      author={Dupont, Johan~L.},
       title={{$K$}-theory obstructions to the existence of vector fields},
        date={1974},
        ISSN={0001-5962,1871-2509},
     journal={Acta Math.},
      volume={133},
       pages={67\ndash 80},
         url={https://doi.org/10.1007/BF02392142},
      review={\MR{425980}},
}

\bib{Fa:time}{article}{
      author={Fathi, Albert},
       title={Time functions revisited},
        date={2015},
        ISSN={0219-8878,1793-6977},
     journal={Int. J. Geom. Methods Mod. Phys.},
      volume={12},
      number={8},
       pages={1560027, 12},
         url={https://doi.org/10.1142/S0219887815600270},
      review={\MR{3400666}},
}

\bib{FaSi}{article}{
      author={Fathi, Albert},
      author={Siconolfi, Antonio},
       title={On smooth time functions},
        date={2012},
        ISSN={0305-0041,1469-8064},
     journal={Math. Proc. Cambridge Philos. Soc.},
      volume={152},
      number={2},
       pages={303\ndash 339},
         url={https://doi.org/10.1017/S0305004111000661},
      review={\MR{2887877}},
}

\bib{GHL:Rie}{book}{
      author={Gallot, Sylvestre},
      author={Hulin, Dominique},
      author={Lafontaine, Jacques},
       title={Riemannian geometry},
     edition={Third},
      series={Universitext},
   publisher={Springer-Verlag, Berlin},
        date={2004},
        ISBN={3-540-20493-8},
         url={https://doi.org/10.1007/978-3-642-18855-8},
}

\bib{G:null}{article}{
      author={Galloway, Gregory~J.},
       title={A note on null distance and causality encoding},
        date={2024},
        ISSN={0264-9381,1361-6382},
     journal={Classical Quantum Gravity},
      volume={41},
      number={1},
       pages={Paper No. 017001, 5},
}

\bib{GH:top}{article}{
      author={Garc\'ia-Heveling, Leonardo},
       title={Topology change with {M}orse functions progress on the
  {B}orde-{S}orkin conjecture},
        date={2023},
        ISSN={1095-0761,1095-0753},
     journal={Adv. Theor. Math. Phys.},
      volume={27},
      number={4},
       pages={1159\ndash 1190},
         url={https://doi.org/10.4310/atmp.2023.v27.n4.a4},
      review={\MR{4761209}},
}

\bib{Gui}{incollection}{
      author={Guilbault, Craig~R.},
       title={Ends, shapes, and boundaries in manifold topology and geometric
  group theory},
        date={2016},
   booktitle={Topology and geometric group theory},
      series={Springer Proc. Math. Stat.},
      volume={184},
   publisher={Springer, [Cham]},
       pages={45\ndash 125},
         url={https://doi.org/10.1007/978-3-319-43674-6_3},
      review={\MR{3598162}},
}

\bib{Haw:time}{article}{
      author={Hawking, S.~W.},
       title={The existence of cosmic time functions},
        date={1969},
        ISSN={00804630},
     journal={Proceedings of the Royal Society of London. Series A,
  Mathematical and Physical Sciences},
      volume={308},
      number={1494},
       pages={433\ndash 435},
         url={http://www.jstor.org/stable/2416157},
}

\bib{HE}{book}{
      author={Hawking, S.~W.},
      author={Ellis, G. F.~R.},
       title={The large scale structure of space-time},
      series={Cambridge Monographs on Mathematical Physics},
   publisher={Cambridge University Press, London-New York},
        date={1973},
      volume={No. 1},
}

\bib{HR}{article}{
      author={Hopf, H.},
      author={Rinow, W.},
       title={Ueber den {B}egriff der vollst\"{a}ndigen
  differentialgeometrischen {F}l\"{a}che},
        date={1931},
        ISSN={0010-2571,1420-8946},
     journal={Comment. Math. Helv.},
      volume={3},
      number={1},
       pages={209\ndash 225},
         url={https://doi.org/10.1007/BF01601813},
      review={\MR{1509435}},
}

\bib{Ho:vf}{article}{
      author={Hopf, Heinz},
       title={Vektorfelder in {$n$}-dimensionalen {M}annigfaltigkeiten},
        date={1927},
        ISSN={0025-5831,1432-1807},
     journal={Math. Ann.},
      volume={96},
      number={1},
       pages={225\ndash 249},
         url={https://doi.org/10.1007/BF01209164},
      review={\MR{1512316}},
}

\bib{Hor}{article}{
      author={Horowitz, Gary~T.},
       title={Topology change in classical and quantum gravity},
        date={1991},
        ISSN={0264-9381,1361-6382},
     journal={Classical Quantum Gravity},
      volume={8},
      number={4},
       pages={587\ndash 601},
         url={http://stacks.iop.org/0264-9381/8/587},
      review={\MR{1100121}},
}

\bib{Hus}{book}{
      author={Husemoller, Dale},
       title={Fibre bundles},
     edition={Third},
      series={Graduate Texts in Mathematics},
   publisher={Springer-Verlag, New York},
        date={1994},
      volume={20},
        ISBN={0-387-94087-1},
         url={https://doi.org/10.1007/978-1-4757-2261-1},
      review={\MR{1249482}},
}

\bib{Jam:spherevfs}{article}{
      author={James, I.~M.},
       title={Whitehead products and vector-fields on spheres},
        date={1957},
        ISSN={0008-1981},
     journal={Proc. Cambridge Philos. Soc.},
      volume={53},
       pages={817\ndash 820},
         url={https://doi.org/10.1017/s0305004100032928},
      review={\MR{105108}},
}

\bib{Ko:vflecture}{book}{
      author={Koschorke, Ulrich},
       title={Vector fields and other vector bundle morphisms---a singularity
  approach},
      series={Lecture Notes in Mathematics},
   publisher={Springer, Berlin},
        date={1981},
      volume={847},
        ISBN={3-540-10572-7},
      review={\MR{611333}},
}

\bib{Lev}{article}{
      author={Levin, V.~L.},
       title={Continuous utility theorem for closed preorders on a metrizable
  {$\sigma $}-compact space},
        date={1983},
        ISSN={0002-3264},
     journal={Dokl. Akad. Nauk SSSR},
      volume={273},
      number={4},
       pages={800\ndash 804},
      review={\MR{728276}},
}

\bib{Mar}{article}{
      author={Markus, L.},
       title={Line element fields and {L}orentz structures on differentiable
  manifolds},
        date={1955},
        ISSN={0003-486X},
     journal={Ann. of Math. (2)},
      volume={62},
       pages={411\ndash 417},
         url={https://doi.org/10.2307/1970071},
      review={\MR{73169}},
}

\bib{Ma72}{article}{
      author={Marsden, Jerrold},
       title={On completeness of homogeneous pseudo-riemannian manifolds},
        date={1972/73},
        ISSN={0022-2518,1943-5258},
     journal={Indiana Univ. Math. J.},
      volume={22},
       pages={1065\ndash 1066},
         url={https://doi.org/10.1512/iumj.1973.22.22089},
}

\bib{Min:LCT}{article}{
      author={Minguzzi, E.},
       title={Lorentzian causality theory},
        date={2019},
     journal={Living Rev. Relativ.},
      volume={22},
      number={3},
}

\bib{Min:cone}{article}{
      author={Minguzzi, Ettore},
       title={Causality theory for closed cone structures with applications},
        date={2019},
        ISSN={0129-055X,1793-6659},
     journal={Rev. Math. Phys.},
      volume={31},
      number={5},
       pages={1930001, 139},
         url={https://doi.org/10.1142/S0129055X19300012},
      review={\MR{3955368}},
}

\bib{ON}{book}{
      author={O'Neill, Barrett},
       title={Semi-{R}iemannian geometry},
      series={Pure and Applied Mathematics},
   publisher={Academic Press, Inc., New York},
        date={1983},
      volume={103},
        ISBN={0-12-526740-1},
        note={With applications to relativity},
}

\bib{Orl:cone}{article}{
      author={Orlitzky, Michael},
       title={Continuity of the conic hull},
        date={2024},
        ISSN={0944-6532,2363-6394},
     journal={J. Convex Anal.},
      volume={31},
      number={1},
       pages={255\ndash 264},
}

\bib{Po:S2}{article}{
      author={Poincar\'e, H.},
       title={Sur les courbes d\'efinies par les \'equations
  diff\'erentielles},
    language={fr},
        date={1885},
     journal={Journal de Math\'ematiques Pures et Appliqu\'ees},
      volume={4e s{\'e}rie, 1},
         url={http://www.numdam.org/item/JMPA_1885_4_1__167_0/},
}

\bib{SaSo}{article}{
      author={Sakovich, A.},
      author={Sormani, C.},
       title={The null distance encodes causality},
        date={2023},
        ISSN={0022-2488},
     journal={J. Math. Phys.},
      volume={64},
      number={1},
       pages={Paper No. 012502, 18},
         url={https://doi.org/10.1063/5.0118979},
      review={\MR{4538267}},
}

\bib{San:gh}{article}{
      author={S\'{a}nchez, Miguel},
       title={Globally hyperbolic spacetimes: slicings, boundaries and
  counterexamples},
        date={2022},
        ISSN={0001-7701,1572-9532},
     journal={Gen. Relativity Gravitation},
      volume={54},
      number={10},
       pages={Paper No. 124, 52},
         url={https://doi.org/10.1007/s10714-022-03002-6},
      review={\MR{4498990}},
}

\bib{SoWo}{article}{
      author={Sorkin, R.~D.},
      author={Woolgar, E.},
       title={A causal order for spacetimes with {$C^0$} {L}orentzian metrics:
  proof of compactness of the space of causal curves},
        date={1996},
        ISSN={0264-9381,1361-6382},
     journal={Classical Quantum Gravity},
      volume={13},
      number={7},
       pages={1971\ndash 1993},
         url={https://doi.org/10.1088/0264-9381/13/7/023},
      review={\MR{1400951}},
}

\bib{SV}{article}{
      author={Sormani, C.},
      author={Vega, C.},
       title={Null distance on a spacetime},
        date={2016},
     journal={Classical Quantum Gravity},
      volume={33},
      number={8},
       pages={085001},
}

\bib{Ste}{book}{
      author={Steenrod, Norman},
       title={The {T}opology of {F}ibre {B}undles},
   publisher={Princeton University Press, Princeton, N. J.},
        date={1951},
      volume={vol. 14.},
      review={\MR{39258}},
}

\bib{Sti:parallel}{article}{
      author={Stiefel, E.},
       title={Richtungsfelder und {F}ernparallelismus in n-dimensionalen
  {M}annigfaltigkeiten},
        date={1935},
        ISSN={0010-2571,1420-8946},
     journal={Comment. Math. Helv.},
      volume={8},
      number={1},
       pages={305\ndash 353},
         url={https://doi.org/10.1007/BF01199559},
      review={\MR{1509530}},
}

\bib{Ma:proc}{proceedings}{
      author={Masushita, Yasuo},
      editor={Suh, Young~Jin},
      editor={Ohnita, Yoshihiro},
      editor={Kim, Byung~Hak},
      editor={Lee, Hyunjin},
       title={Lecture notes on indefinite manifolds \& related topics},
      series={Proceedings of the 21st {I}nternational {W}orkshop on {H}ermitian
  {S}ymmetric {S}paces and {S}ubmanifolds and 14th {RIRCM}-{OCAMI} {J}oint
  {D}ifferential {G}eometry {W}orkshop},
        date={2017},
        note={Held at Kyungpook National University, Daegu Republic of Korea,
  October 11--13, 2017},
}

\bib{Tho:vfs}{article}{
      author={Thomas, Emery},
       title={Vector fields on manifolds},
        date={1969},
        ISSN={0002-9904},
     journal={Bull. Amer. Math. Soc.},
      volume={75},
       pages={643\ndash 683},
         url={https://doi.org/10.1090/S0002-9904-1969-12240-8},
      review={\MR{242189}},
}

\bib{Ve}{misc}{
      author={Vega, Carlos},
       title={Spacetime distances: an exploration},
   publisher={arXiv:2103.01191},
        date={2021},
}

\bib{Whi:top}{article}{
      author={Whitney, Hassler},
       title={Topological properties of differentiable manifolds},
        date={1937},
        ISSN={0002-9904},
     journal={Bull. Amer. Math. Soc.},
      volume={43},
      number={12},
       pages={785\ndash 805},
         url={https://doi.org/10.1090/S0002-9904-1937-06642-0},
      review={\MR{1563640}},
}

\end{biblist}
\end{bibdiv}

\end{document}